\documentclass[10pt,a4,twoside,rm]{article}
\usepackage[left=4.5cm,top=4cm,right=4.5cm,bottom=2cm]{geometry}
\usepackage{amsmath, amsthm, amssymb}
\usepackage{graphicx}
\usepackage[all, knot]{xy}
\usepackage{amssymb}
\usepackage{fancyhdr}
\usepackage{titlesec}
\usepackage{multirow}
\newtheorem{Thm}{Theorem}

\newtheorem{lemma}{Lemma}

\newtheorem{defi}[Thm]{Definition}
\newtheorem{conj}{Conjecture}
\newtheorem*{conj*}{Conjecture}

\newcommand{\Si}{\mathfrak{S}}

\newcommand{\Stab}{\operatorname{Stab}}
\newcommand{\Par}{\operatorname{Par}}

\newcommand{\Hom}{\operatorname{Hom}}

\newcommand{\Span}{\operatorname{span}}

\newcommand{\Ind}{{\cal I}}

\newcommand{\rad}{\operatorname{rad}}

\newcommand{\Std}{\operatorname{Std}}

\newcommand{\C}{{\mathbb C}}
\newcommand{\Z}{{\mathbb Z}}
\newcommand{\Q}{{\mathbb Q}}
\newcommand{\F}{{\mathbb F}_p}

\newcommand{\LL}{{\cal L}}
\newenvironment{pf*}{\noindent \textit{Proof:} }{\quad \hfill $\square$ \newline}

\edef\savecatcodeat{\the\catcode`@}
\catcode`\@=11

\def\tb@ifSpecChars#1#2{#1}
\def\tb@ifNoSpecChars#1#2{#2}

\def\tableau{%
  \bgroup
  \@ifstar{\let\Tif\tb@ifNoSpecChars\tb@tableauB}
          {\let\Tif\tb@ifSpecChars\tb@tableauB}}

\def\tb@tableauB{
  \@ifnextchar[{\tb@tableauC}{\tb@tableauC[]}}

\def\tb@tableauC[#1]{\hbox\bgroup%
    \let\\=\cr
    \def\bl{\global\let\tbcellF\tb@cellNF}%
    \def\tf{\global\let\tbcellF\tb@cellH}
    \dimen2=\ht\strutbox \advance\dimen2 by\dp\strutbox%
    \ifx\baselinestretch\undefined\relax%
    \else%
       \dimen0=100sp \dimen0=\baselinestretch\dimen0%
       \dimen2=100\dimen2 \divide\dimen2 by\dimen0%
    \fi%
    \let\tpos\tb@vcenter
    \tb@initYoung
    \tb@options#1\eoo
    \let\arrow\tb@arrow%
    \dimen0=\Tscale\dimen2%
    \dimen1=\dimen0 \advance\dimen1 by \tb@fframe%
    \lineskip=0pt\baselineskip=0pt
    \def\tb@nothing{}%
    \def\endcellno{$\rss\egroup\bss\egroup}
    \def\endcell{\endcellno\kern-\dimen0}
    \def\begincell{\vbox to\dimen0\bgroup\vss\hbox to\dimen0\bgroup\hss$}%
    \let\overlay\tb@overlay%
    \let\fl\tb@fl%
    \let\lss\hss\let\rss\hss\let\tss\vss\let\bss\vss
    \def\mkcell##1{
        \let\tbcellF\tb@cellD
        \def\tb@cellarg{##1}
        \ifx\tb@cellarg\tb@nothing\let\tb@cellarg\tb@cellE\fi%
        \begincell\tb@cellarg\endcellno
        \tbcellF}
    \let\savecellF\tbcellF
     \Tif{\catcode`,=4\catcode`|=\active}{}\tb@tableauD}%

\let\tb@savetableauD\tableauD
{
    \catcode`|=\active \catcode`*=\active \catcode`~=\active%
    \catcode`@=\active
\gdef\tableauD#1{%
  \Tif{
    \mathcode`|="8000 \mathcode`*="8000%
    \mathcode`~="8000 \mathcode`@="8000%
    \def@{\bullet}%
    \let|\cr
    \let*\tf
    \let~\sk
  }{}%
  \tpos{\tabskip=0pt\halign{&\mkcell{##}\cr#1\crcr}}%
  \global\let\tbcellF\savecellF
  \egroup
  \egroup}
}
\let\tb@tableauD\tableauD
\let\tableauD\tb@savetableauD
\let\tb@savetableauD\undefined

\def\tb@options#1{\ifx#1\eoo\relax\else\tb@option#1\expandafter\tb@options\fi}

\def\tb@option#1{%
  \if#1t\let\tpos\tb@vtop\fi
  \if#1c\let\tpos\tb@vcenter\fi
  \if#1b\let\tpos\vbox\fi
  \if#1F\tb@initFerrers\fi
  \if#1Y\tb@initYoung\fi
  \if#1s\tb@initSmall\fi
  \if#1m\tb@initMedium\fi
  \if#1l\tb@initLarge\fi
  \if#1p\tb@initPartition\fi
  \if#1a\tb@initArrow\fi
}

\def\tb@vcenter#1{\ifmmode\vcenter{#1}\else$\vcenter{#1}$\fi}

\def\tb@vtop#1{\hbox{\raise\ht\strutbox\hbox{\lower\dimen0\vtop{#1}}}}

\def\tb@initPartition{\def\Tscale{.3}}
\def\tb@initSmall{\def\Tscale{1}}
\def\tb@initMedium{\def\Tscale{2}}
\def\tb@initLarge{\def\Tscale{3}}

\def\tb@initArrow{\dimen2=1.25em}

\def\tb@initYoung{%
  \def\tb@cellE{}
  \let\tb@cellD\tb@cellN
  \def\sk{\global\let\tbcellF\tb@cellNF}}
\def\tb@initFerrers{%
  \def\tb@cellE{\bullet}
  \let\tb@cellD\tb@cellNF
  \def\sk{\bullet}}

\tb@initMedium

\def\tb@sframe#1{%
  \vbox to0pt{
    \vss
    \hbox to0pt{%
      \hss
      \vbox to\dimen1{
        \hrule depth #1 height0pt
        \vss
        \hbox to\dimen1{
          \vrule width #1 height\dimen1
          \hss
          \vrule width #1
          }%
        \vss
        \hrule height #1 depth 0in
        }%
      \kern-\tb@hframe
      }%
    \kern-\tb@hframe}}

\def\tb@hframe{.2pt}\def\tb@fframe{.4pt}\def\tb@bframe{1.2pt}
\def\tb@cellH{\tb@sframe{\tb@bframe}}       
\def\tb@cellNF{}                            
\def\tb@cellN{\tb@sframe{\tb@fframe}}       
\let\tbcellF\tb@cellN                       

\def\tb@rpad{1pt}
\def\tb@lpad{1pt}
\def\tb@tpad{1.8pt}
\def\tb@bpad{1.8pt}

\def\tb@overlay{\endcell\@ifnextchar[{\tb@overlaya}{\begincell}}
\def\tb@overlaya[#1]{\vbox to\dimen0\bgroup%
  \tb@overlayoptions#1\eoo%
  \tss\hbox to\dimen0\bgroup\lss}
\def\tb@overlayoptions#1{\ifx#1\eoo\relax\else\tb@overlayoption#1\expandafter\tb@overlayoptions\fi}

\def\tb@overlayoption#1{
  \if#1t\def\tss{\vskip\tb@tpad}\let\bss\vss\fi
  \if#1c\let\tss\vss\let\bss\vss\fi
  \if#1b\def\bss{\vskip\tb@bpad}\let\tss\vss\fi
  \if#1l\def\lss{\hskip\tb@lpad}\let\rss\hss\fi
  \if#1m\let\lss\hss\let\rss\hss\fi
  \if#1r\def\rss{\hskip\tb@rpad}\let\lss\hss\fi
}

\def\tb@fl{\endcell\begincell\vrule depth 0pt width \dimen0 height \dimen0 \endcell\begincell}

\@ifundefined{diagram}{}{
\def\tb@arrowpad{.5}

\newoptcommand{\tb@arrow}{\@ne}[2]{%
  \endcell
   \begingroup%
   \let\dg@getnodesize\tb@getnodesize
   \dg@USERSIZE=#1\relax%
   \ifnum\dg@USERSIZE<\@ne \dg@USERSIZE=\@ne \fi%
   \dg@parse{#2}%
   \dg@label{\tb@draw{#1}{#2}}}

\def\tb@getnodesize#1#2#3#4#5{\dimen3=\tb@arrowpad\dimen2 #4=\dimen3 #5=\dimen3\relax}
\def\tb@getnodesize#1#2#3#4#5{\ifnum#2=0\ifnum#3=0\tb@getnodesizetail{#4}{#5}\else\tb@getnodesizehead{#4}{#5}\fi\else\tb@getnodesizehead{#4}{#5}\fi}
\def\tb@getnodesizetail#1#2{\dimen3=.5\dimen2 #1=\dimen3 #2=\dimen3}
\def\tb@getnodesizehead#1#2{\dimen3=.5\dimen2 #1=\dimen3 #2=\dimen3}

\def\tb@draw#1#2#3#4{%
        \dg@X=0\dg@Y=0\dg@XGRID=1\dg@YGRID=1\unitlength=.001\dimen0%
        \dg@LBLOFF=\dgLABELOFFSET \divide\dg@LBLOFF\unitlength%
        \dg@drawcalc
        \begincell
        \let\lams@arrow\tb@lams@arrow
        \begin{picture}(0,0)\begingroup\dg@draw{#1}{#2}{#3}{#4}\end{picture}%
        \endcell
        \endgroup
        \begincell}
}

\def\tb@lams@arrow#1#2{%
 \lams@firstx\z@\lams@firsty\z@
 \lams@lastx#1\relax\lams@lasty#2\relax
 \lams@center\z@
 \N@false\E@false\H@false\V@false
 \ifdim\lams@lastx>\z@\E@true\fi
 \ifdim\lams@lastx=\z@\V@true\fi
 \ifdim\lams@lasty>\z@\N@true\fi
 \ifdim\lams@lasty=\z@\H@true\fi
 \NESW@false
 \ifN@\ifE@\NESW@true\fi\else\ifE@\else\NESW@true\fi\fi
 \ifH@\else\ifV@\else
  \lams@slope
  \ifnum\lams@tani>\lams@tanii
   \lams@ht\ten@\p@\lams@wd\ten@\p@
   \multiply\lams@wd\lams@tanii\divide\lams@wd\lams@tani
  \else
   \lams@wd\ten@\p@\lams@ht\ten@\p@
   \divide\lams@ht\lams@tanii\multiply\lams@ht\lams@tani
  \fi
 \fi\fi
 \ifH@  \lams@harrow
 \else\ifV@ \lams@varrow
 \else \lams@darrow
 \fi\fi
}

\catcode`\@=\savecatcodeat
\let\savecatcodeat\undefined

\titleformat{\section}[hang]{\sc}{\thesection.}{0.5cm}{\filcenter}
\titleformat{\subsection}[hang]{\it}{\thesubsection.}{0.2cm}{\filright}

\begin{document}

\title{\bf \normalsize Projective modules for the symmetric group and Young's seminormal form.}
\author{\sc  steen ryom-hansen\thanks{Supported in part by FONDECYT grant 1121129, 
 by Programa Reticulados y Simetr\'ia
and by the MathAmSud project OPECSHA 01-math-10} }
\date{}   \maketitle
\begin{abstract}
We study the representation theory of the symmetric group $ \Si_n $ in positive characteristic $p$.
Using features of the LLT-algorithm we give a
conjectural description of the projective cover $ P(\lambda) $ of 
the simple module $ D(\lambda) $ where $ \lambda $ is a $p$-restricted partition
such that all ladders of the corresponding ladder partition are of order less than $p$.
Inspired by the recent theory of KLR-algebras, we explain an 
algorithm that allows us to verify this conjectural description for $ n \leq 15 $, at least. 
\end{abstract}

\section{Introduction} 
Computing the decomposition matrices of the symmetric groups is one of the big open problems in representation theory. 
Until recently, the work in this area was guided by the James' conjecture
which says that the decomposition matrices should coincide with those for the Hecke algebras of type $A$ at roots of unity, 
for the prime in a certain range.
In particular, the conjecture predicts the decomposition numbers to be given by certain Kazhdan-Lusztig polynomials.

\medskip
Two of the recent major developments in representation theory, the Brundan and Kleshchev's isomorphism Theorem  
together with Elias and Williamson's algebraic proof of Soergel's conjecture, could have
be seen as evidence in favor of the James' conjecture. But even more recently, Williamson changed the 
subject dramatically by publishing a preprint that gave counterexamples to the James' conjecture.

\medskip
In this paper we combine the classical theory of Young's seminormal form, the Lascoux, Leclerc and Thibon algorithm and
Brundan and Kleshchev's isomorphism Theorem 
to make the $ \Si_n $-representation theory look formally
like the representation theory of Soergel bimodules. This opens up a new perspective on the representation theory of 
$ \Si_n $ which we believe will be fruitful in the future.

\medskip
Let us explain the contents of the paper in more detail. Let $p $ be a prime and let 
$ \F $ be the field of $p$ elements. Denote by $\Par_{res,n}$ 
the set of $ p$-restricted partitions of $ n$. As is well known, it parametrizes the simple 
$ \F \Si_n $-modules so that we for $ \lambda \in \Par_{res,n} $ have
a simple $ \F \Si_n $-module $ D(\lambda)$ together with its projective 
cover $P(\lambda) $.

\medskip
After section 2, which is devoted to setting up the notation, we construct in section 3
an  
idempotent $ \widetilde{e}_{\lambda} \in \F \Si_n$ that plays an important role 
throughout the paper. 
The main ingredient for $ \widetilde{e}_{\lambda} $ is 
Murphy's tableau class idempotent for the ladder class of $ \lambda $.

\medskip
Defining $ \widetilde{A(\lambda)}:= \F \Si_n \widetilde{e}_{\lambda}$ we obtain a projective 
$\F \Si_n$-module, but it is decomposable in general, that is $ \widetilde{A(\lambda)} \not= P(\lambda)$.
On the other hand, we show in Theorem 
{\ref{expansion}}  
that there is triangular expansion of the form
\begin{equation}{\label{firsttriangular}} \widetilde{A(\lambda)} = P(\lambda) \oplus \bigoplus_{\mu, \mu \rhd \lambda } 
\, P(\mu)^{\oplus m_{\lambda \mu}}
\end{equation}
for certain nonnegative integers $ m_{\lambda \mu} $ 
where $ \rhd $ is the usual dominance order on partitions.

\medskip
In section 4 we consider
the Lascoux, Leclerc and Thibon (LLT) algorithm which gives a way of calculating the global crystal 
basis $ \{ G(\lambda) \,  | \,  \lambda \in \Par_{res,n} \} $, 
for the basic submodule $ {\cal M}_q $ of the $q$-Fock space. 
An important ingredient for this algorithm is given by 
certain elements $ A(\lambda) $ of $ {\cal M}_q $
that LLT called `the first approximation of the global basis'. 
In fact, their algorithm is a triangular recursion based on these elements.
We use this to observe that 
they satisfy the following triangular expansion property
$$ {A(\lambda)} = G(\lambda) + \sum_{\mu, \mu \rhd \lambda } 
\,n_{\lambda \mu}(q) G(\mu)  $$
where $ n_{\lambda \mu}(q) \in \Z[q,q^{-1}]$.

\medskip
Our main point is now to consider $ A(\lambda) $ as an object of interest in itself, and not just a tool for 
calculating $ G(\lambda)$. In this spirit
we conjecture that ${A(\lambda)}$ should be categorified by $ \widetilde{A(\lambda)} $, 
or to be more precise that we should have
\begin{equation}{\label{mainconjectureintro}}
 n_{\lambda \mu}(1)= m_{\lambda \mu}.
\end{equation}
This formula is the main theme of our paper.
In Theorem {\ref{mainJames}} we show 
that it implies 
James' conjecture.

\medskip
In section 5 we describe a method for verifying ({\ref{mainconjectureintro}})
for $ n $ not too big. On the other hand, in 
view of Theorem {\ref{mainJames}}
and 
Williamson's counterexamples, ({\ref{mainconjectureintro}}) cannot be true for all $ p  $ such that $ n < p^2 $, that is within the range 
for James' conjecture. Williamson's smallest counterexample is big, $ n =467874  $ and $ p= 839$, 
and so we do not speculate on the true range of validity for (\ref{mainconjectureintro}). On the other hand, 
our method for verifying (\ref{mainconjectureintro}) is of interest in itself. It is based
on Brundan and Kleshchev's isomorphism between $ \F \Si_n $ and ${ \cal  R}_n $, the  
cyclotomic Khovanov, Lauda and Rouquier (KLR) algebra of type $A$.
Here $ \widetilde{e}_{\lambda} $ is  
closely related to 
the KLR-idempotents and so $ \widetilde{e}_{\lambda } S(\mu)$ 
identifies with the symmetrized generalized eigenspace for the action of the Jucys-Murphy elements in the Specht module $ S(\mu)$.
We must calculate the $p$-rank of the canonical bilinear form 
$ \langle \cdot, \cdot \rangle_{\mu} $ on $ S(\mu)$ on the restriction to $ \tilde{e}_{\lambda } S(\mu)$.
To do this we rely on our results from [RH3] on 
the compatibility of the `intertwining elements' from Brundan and Kleshchev's work with 
Young's seminormal form. They allow us to describe the action of the KLR-generators $ \psi_i $ 
completely in terms of Young's seminormal form.
Our partial verification of (\ref{mainconjectureintro}) follows from this.

\medskip
In the final section 6 of the paper, we take the relationship with the KLR-algebra one step 
further. Indeed, one of the important aspects of the KLR-algebra 
is the fact that it is 
a $ \Z$-graded algebra in a nontrivial way and hence it is possible to speak of graded modules
over it. 
By comparison with certain idempotents that occur naturally in the nilHecke algebra, 
we show in Theorem {\ref{idempotents}} that the idempotent $ \tilde{e}_{\lambda } $
is a homogeneous idempotent of $ { \cal  R}_n $. In particular, $ \widetilde{A(\lambda)} $
admits a grading and from this it follows from general theory 
that also $ P(\lambda) $ admits a grading. 

\medskip
It is a pleasure to thank the referee for his/hers useful comments.


\section{Basic Notation and a couple of Lemmas}
Let $ p > 2 $ be a prime and let $ R $ be the localization of $ \Z $ at $ p$.
Let $ \Si_n $ be the symmetric group on $ n $ letters and 
write $ \sigma_i := (i-1,i) $. 
We are interested in the representation theory 
of $ \Si_n $ over the finite field $ \F =R/pR $.

\medskip 
Over $\Q $, the irreducible representations of $\Si_n $ are parametrized by 
the set $ \Par_{n}$ of partitions of $n$, that is the set of 
nonincreasing sequences of positive integers 
$ \lambda =(\lambda_1, \ldots, \lambda_k ) $ with sum $n$. Over $ \F $ they are
parametrized by the set of 
$p$-restricted partitions $ \Par_{res,n} $, consisting of those $ \lambda \in \Par_n $ 
that satisfy $ \lambda_{i} - \lambda_{i+1} < p $ for all $ i $ where by convention $ \lambda_i = 0 $ for $ i \geq k+1$. 
For $ \lambda \in \Par_n $ we denote by $ S(\lambda) $ the Specht module for $ R  \Si_n$, 
see below for the precise definition.
In general, for an $ R \Si_n$-module $ M $ 
we denote by $ \overline{M} := M \otimes_R \F $ the $ \F \Si_n$-module
obtained by reduction modulo $p$, but sometimes, when there is no risk of confusion, we also 
refer to it simply as $ M $.
There is a 
bilinear, symmetric $ \Si_n$-invariant form $  \langle \cdot, \cdot \rangle_{\lambda} $ on 
$ S(\lambda) $ 
which is nonzero iff $ \lambda \in \Par_{res,n}$ and we obtain the parametrization of the simple 
modules for $ \F \Si_n$ via $ \lambda \in \Par_{res,n} \mapsto D(\lambda) 
:=\overline{S(\lambda)}/ \rad \langle \cdot, \cdot \rangle_{\lambda} $.

\medskip
In the paper we shall be specially interested in the projective covers of the simple modules.
For $ \lambda \in \Par_{res,n}$ 
we denote by $ P(\lambda ) $ the projective cover of $ D(\lambda) $.
By definition, $ P(\lambda) $ is the unique indecomposable projective $ \F \Si_n$-module such that 
$ D(\lambda) $ is a homomorphic image of $ P(\lambda) $.
By general theory, $ P(\lambda) $ is of the form $ P(\lambda) = \F \Si_n e_{\lambda} $ for some 
idempotent $ e_{\lambda} \in \F \Si_n$. Unfortunately, there is in general no concrete 
description of $ e_{\lambda}$.

\medskip
A 
partition $ \lambda =  (\lambda_1, \ldots, \lambda_k ) \in \Par_n$ 
is represented graphically via its Young diagram. It consists of 
$ k $, left aligned, files of boxes, called nodes, in the plane,
with the first file containing $ \lambda_1 $ nodes, the second file containing 
$ \lambda_2 $ nodes and so on. The nodes are indexed using matrix convention, with the $[i,j]$'th node 
situated in the $j$'th column of the $i$'th file.
For $ \lambda \in \Par_n $, a $ \lambda$-tableau $ t $ is a filling 
of the nodes of $ \lambda $ with the numbers $ \{1, 2, \ldots, n \} $.
We write 
$ t[i,j] = k $ if the $ [i,j] $'th node of $ t $ is filled with $ k$
and $ c_{ t}( k ) = j-i $ if $ t[i, j ] = k$. Then $ c_{ t}( k )$ is the content of $ t $ 
at $k$, whereas its image in $ \F $, denoted $ r_{ t}( k )$, is the $p$-residue of $t $ at $ k$.
For $ k \in \{1, 2, \ldots, n \} $ 
we define $ t(k) := [i,j] $ where $ t[i,j] = k $. 
A tableau $ t $ is called standard if $ t[i, j ] \leq t[i, j+1 ] $ and $ t[i, j ] \leq t[i+1, j ] $
for all relevant $ i, j $.
The set of standard tableaux of partitions of $ n$ is denoted $ \Std(n) $ and the set of 
standard tableaux with underlying partition $ \lambda $ is denoted $ \Std(\lambda)$. 
For $ \lambda \in \Par_n$ and $ t \in  \Std(\lambda)$-tableau we 
write $ Shape(t) := \lambda $.

\medskip
Let $ t $ be a $ \lambda $-tableau with node $ [i,j] $. 
The $ [i, j]$-hook consists 
of the nodes of the Young diagram of $ \lambda $ situated to the right and below the $ [i, j ] $ node and its 
cardinality is called the hook-length $ h_{ij} $. The product of all hook-lengths 
is denoted $ h_{\lambda} $.
The hook-quotient of the tableau $ t \in \Std(\lambda) $ at $ n $ is the number 
$ \gamma_{tn} = \prod  \frac{h_{ij}}{h_{ij}-1} $ 
where the product is taken over 
all nodes in the row of $ \lambda $ that contains $ n $, omitting hooks of length one. 
For a general $ i $, we define $ \gamma_{ti} $ similarly, by first deleting from $ t $ the nodes containing 
$ i+1, i+2, \ldots , n $. 
Finally we define $ \gamma_t = \prod_{i=2}^n \gamma_{ti} $.

\medskip
Let us recall the
combinatorial concepts of ladders and ladder tableaux that play an important role for the
LLT-algorithm, although we shall use conventions that are dual to the ones of [LLT].
Let $ \mu $ be a $p$-restricted partition. The 'ladders' of $\mu$ are the
straight 'line segments' through the Young diagram of $ \mu $ with 'slope' $ 1/(p-1) $, 
that is the subsets of the nodes of $ \mu $ of the form $ { \cal L}_b :=\{ \, [i,j] \, | \, j= b-(p-1)(i-1) \} $.
If $\mu \in \Par_{res,n}$ we have that the ladders are 'unbroken', that is $ \pi_1( {\cal L}_a) $ is 
of the form $\{ q ,q+1, a +2, \ldots, r \} $ for some $ q< r $ where $ \pi_1 $ is the first projection.
We say that $ { \cal L}_b $ is smaller than $ { \cal L}_{b_1} $ if $ b < b_1$. The ladder tableau 
$ \mu_{lad} $ of $ \mu$ is defined as the $\mu$-tableau with the numbers $ 1, 2, \ldots, n $ filled in one
ladder at the time, 
starting with the smallest ladder and continuing successively upwards, 
the numbers being filled in from top to bottom in each ladder.
Note that the residues are constant on each ladder.

\medskip
For a partition $ \mu $, the $p$-residue diagram $ res_{\mu} $ is obtained by writing the residue
$ r_{ t}( k )$ 
in the $[i,j]$'th node of the Young diagram of $ \mu$. For example, 
if $ \mu= (6, 5, 3,1) $ and $ p= 3 $ then the residue diagram and ladder tableau are as follows
$$ res_{\mu}= 
{\small{\tableau[scY]{0,1, 2, 0, 1, 2| 2, 0, 1 , 2,0 | 1, 2 , 0| 0| }}}, \, \, \, \, \, \, \, \, \, \, \, \, \, \, \, \, \, \, 
\mu_{lad} = {\small{\tableau[scY]{1,2, 3, 5, 7, 10| 4, 6,8, 11 , 13 | 9, 12 , 14|15| }}}
$$
\medskip \noindent
with ladders $ {\cal L}_1 = \{1\}, {\cal L}_2= \{2\}, {\cal L}_3= \{3,4\}, 
{\cal L}_4 = \{5,6\}, 
{\cal L}_5 = \{7,8,9\}, {\cal L}_6 = \{10,11,12\} $ and $ {\cal L}_7 = \{13,14,15\} $. 
We denote by $ {\bf i}_{ \,lad, \mu } $ 
the residue sequence given by the ladder tableau for $ \mu $. In the above example
it is 
$$ {\bf {i}}_{\, lad,\mu} = (0,1,2,2, 0,0, 1,1,1,2,2,2,0,0,0 ). $$
The ladders define a sequence of subpartitions $ \mu_{ lad,\, \le 1 }, \ldots,  
\mu_{ lad, \, \le m	 } $ of $ \mu $ where $ \mu_{ lad, \, \le k   } $ is defined as
the union of the ladders 
$ {\cal L}_1, {\cal L}_2, \ldots, {\cal L}_k $.

\medskip
We define positive integers $ n_0, \ldots, n_m $ by  $ n_0 := 0 $ and 
\begin{equation}{\label{ladder-limits}}
n_k := | {\cal L}_1 | + | {\cal L}_2 | + \ldots + | {\cal L}_k|.
\end{equation}
We may then 
introduce the ladder group $ \Si_{lad, \mu}  \leq \Si_n$ for $ \mu $ as $ \Si_{lad, \mu}:= \prod_k \Si_{ {\cal L}_{k} }, $
where $\Si_{ {\cal L}_{k} } $ is the symmetric group on the letters  
$ n_{k-1}+1, \ldots, n_{k} $.

\medskip

The dominance order $ \unlhd $ on partitions is defined by 
$$  \lambda   \unlhd \mu   \mbox{ if }
\sum_{i=1}^{m} \lambda_i \leq \sum_{i=1}^{m} \mu_i \, \, \mbox{   for } m = 1, 2, \ldots ,
\mbox{min}(k, l) $$
for $ \lambda = (\lambda_1, \ldots, \lambda_k) $ and 
$ \mu = (\mu_1, \ldots, \mu_l) $. 
When $ \lambda $ is used as a subscript where a tableau is expected, it refers 
to the unique maximal $\lambda$-tableau $ t^{\lambda}$, having the numbers $ \{ 1,\ldots, n \} $ 
filled in along the rows.
The dominance order extends to tableaux by considering them as series of partitions. 

\medskip
Following Murphy in [Mu83], we define an equivalence relation on the set of all standard tableaux via 
$ t \sim_p s $ if $ r_t(k) = r_s(k) \mbox{ mod} \, p $ for all $ k$.
The classes of $ \sim_p $ are called 
tableaux classes. The tableau class containing $ t $ is 
denoted $[t] $. 
The tableaux classes are given by residue sequences, that is elements of $ (\F)^n$, 
although a given residue sequence $ {\bf i } \in (\F)^n $ may give rise to the empty class. 
The ladder tableaux 
are 'minimal' in their classes in the sense of the following Lemma.
Note that throughout we use the convention that 
$ \Si_n $ acts on the left on tableaux by place permutations.
\begin{lemma}{\label{ladder}}
Assume that $ \lambda $ is $ p$-restricted. Then if $ t \in [  \lambda_{lad} ]$ we have that 
either $ Shape(t) \rhd \lambda$ 
or $ Shape(t) = \lambda $ and $ t = \sigma \lambda_{lad}  $ for $ \sigma \in \Si_{lad,\lambda}$.
\end{lemma}
\begin{pf*}
Omitted.
\end{pf*}

Suppose that $\mu \in \Par_n $. A node of $ \mu $ is called removable if it can be removed from 
$ \mu$ with the result being the diagram of a partition $\lambda$. Dually, 
that node is called an addable 
node of $ \lambda$. It is called an $i$-node if its $p$-residue is $i$.

\begin{lemma}{\label{laddergroupaction}}
Assume that $ \mu $ is $ p$-restricted
and that $$ T_{\mu \lambda} :=\{ t \in [\mu_{lad}] \, | \, Shape(t) = \lambda \} \neq \emptyset.$$
Then the ladder group 
$ \Si_{lad, \mu} $ acts faithfully on $ T_{\mu \lambda} $.
\end{lemma}
\begin{pf*}
\medskip
In general, 
we may think of the tableau class $ [ t] $ in an 
algorithmic  way. Indeed, setting $ {\bf i}^t : =
i_1 i_2 \ldots i_n \in (\F)^n$ where $ i_k := r_{ t}( k )$ 
we obtain the tableaux in $ [ t] $ by starting with the one-node partition, 
to which we add in all possible ways an addable $ i_2 $-node. 
For each arising partition, we add in all possible ways an addable $ i_3 $-node and so on. The 
set of tableaux that arises in this way after $n$ steps is exactly $[t]$.
From this, it is clear that 
$ \Si_{lad, \mu} $ acts faithfully on $ T_{\mu \lambda} $.
\end{pf*}

\medskip
For $ k = 1,2, \ldots, n $ the Jucys-Murphy elements $ L_k \in {\mathbb	 Z} \Si_n $ are
defined by
$$ L_k := (1, k ) + (2,k ) + \ldots + (k-1, k) $$
with the convention that $ L_1  := 0 $. 
An important application of the $L_k $ is the construction of 
orthogonal idempotents 
$ E_t \in { \mathbb Q} \Si_n $, the Jucys-Murphy idempotents, indexed by tableaux $ t $, 
that can be used to derive Young's seminormal form. 
Their construction is as follows 
$$ E_t := \prod_{  \{c \, |    -n \, < \,  c \,  < n \} }  \prod_{  \{ \, i \, |  c_{ t}(i) \not= c  \} \,} 
\frac{ L_i -c }{ c_{t}(i) -c}.$$
For $ t $ standard we have $ E_t  \not= 0 $, whereas 
for $ t $ nonstandard either $ E_t = 0 $, or 
$ E_t = E_s $ for some standard tableau $ s $ related to $ t$.
Running over all standard tableaux, the $ E_t $ form a set of primitive and 
complete idempotents for $ \Q \Si_n $, that is their sum is $1$.
Moreover, they are eigenvectors for the action of the Jucys-Murphy operators in $ {\mathbb Q} \Si_n $,
since 
\begin{equation}{\label{Mu_page_506}}
(L_k-c_t(k)) E_t =0 \mbox{ or equivalently }
L_k = \sum_{ t \in \Std(n)} c_t(k) E_t. \end{equation}

\medskip
For $ \lambda \in  \Par_n   $, we 
let $ \Stab_{\lambda} $ denote the row stabilizer of $ t^{\lambda} $ and define $ x_{\lambda} $ and 
$ y_{\lambda} $ as the following elements of $ R\Si_n$ 
$$ x_{\lambda} = \sum_{ \sigma \in \Stab_{\lambda} }   \sigma \,\,\,\, \mbox{and} \, \, \, \, 
y_{\lambda} = \sum_{ \sigma \in \Stab_{\lambda} } (-1)^{| \sigma | } \sigma $$
where $ | \sigma | $ is the sign of $ \sigma $. 
For $ t \in \Std(\lambda) $, we define the associated element $ d(t) \in \Si_n $ by
$$ d(t) t^{\lambda}  = t.$$
Then for pairs of standard $ (s, t) $ of $ \lambda$-tableaux, Murphy's standard basis and dual standard basis, mentioned above, 
consist of 
the elements
$$ x_{s t } = d(s) x_{\lambda } d(t)^{-1}  \,\,\,\, \mbox{and} \, \, 
y_{s t } = d(s) y_{\lambda } d(t)^{-1}. $$
They are bases for $ R \Si_n$ and also for $ \F \Si_n$. 
Set $$ (R \Si_n)^{ > \lambda} := \Span_{R} \{ x_{st} | Shape(s) > \lambda \}. $$
Then $ (R \Si_n)^{ > \lambda} $ is an ideal of $ R\Si_n $ and the Specht module $ S(\lambda)$, mentioned 
above, is the span of 
$ \{ x_{s\lambda} + (R \Si_n)^{ > \lambda} \,| \,s \in \Std(\lambda) \,\}$. 
These elements form an $ R$-basis for $ S(\lambda)$.
Define $ S(\lambda)_{\Q} := S(\lambda) \otimes_{\Z} \Q$ and let 
$ \xi_{st} := E_s x_{st} E_t $. Then 
$ \{\, \xi_{st}  \,|\, (s,t) \in \Std(\lambda)^2, \,\lambda \in \Par_n\, \}$
is the seminormal basis for $ \Q \Si_n$ and, moreover,
$ \{ \xi_{s\lambda} | s \in \Std(\lambda) \}$ is a basis for $ S(\lambda)_{\Q}$.

\medskip
The action of $ \Si_n $ on the standard basis $ \{ x_{s \lambda} \} $ 
is given by a recursion using the Garnir relations, whereas the action of $ \Si_n $ on 
$ \{ \xi_{s \lambda} \} $ 
is given by the following formulas, that appear for example in Theorem 6.4 of [Mu93]
(in the more general context of Hecke algebras, but note the
sign error there: the expression for $ h $ should be replaced by $ -h$).
\begin{Thm}{\label{Murphyysf}}
Let $ h=  c_{ s}( i-1 ) -c_{ s}( i )  $ be the radial distance between the $i-1$ and $ i$-nodes of 
$ s \in \Std(\lambda)$. Let $ t := \sigma_i  s  $ where still $ \sigma_i = (i-1,i)$.
Then the action of $ \sigma_i  $ on $ \xi_{s \lambda} $ is given by the formulas
\begin{equation}{\label{YSFMurphy}}
\sigma_i \xi_{s \lambda}:= \left\{ 
\begin{array}{ll}
\xi_{s \lambda} & \mbox{ if } h = -1 \, \,\,  (i-1 \mbox{ and } i \mbox{ are in same row})  \\
-\xi_{s \lambda} & \mbox{ if } h = 1 \, \,\, (i-1 \mbox{ and } i \mbox{ are in same column}) \\
-\frac{1}{h} \, \xi_{s \lambda} + \xi_{t \lambda}   & \mbox{ if } h > 1 \, \,\, (i-1 \mbox{ is above } i)  \\
-\frac{1}{h} \, \xi_{s \lambda} +\frac{h^2 -1}{h^2} \, \xi_{t \lambda}  & \mbox{ if } h < -1 \, \,\, (i-1 \mbox{ is below } i).
\end{array}
\right.
\end{equation}
\end{Thm}

\section{Projective modules}
For the results of this section,  we first need a description of Robinson's $i$-induction functor 
in terms of the Jucys-Murphy idempotents. There are related descriptions 
available in the literature, see for example [HuMa3], but our description of 
the 'divided power' functor seems to be new. It relies on a result from our recent paper [RH3]. 

\medskip
Recall that the blocks for $  {\mathbb  F}_p \Si_n  $ are given by the 
Nakayama conjecture (which is a Theorem). Murphy showed in [Mu83] how to describe the 
corresponding block idempotents in terms of the Jucys-Murphy idempotents $E_t$. 
Indeed,
let $T=[t]$ be the class of 
$ t \in \Std(n)$ under $ \sim_p$ and consider for $ \lambda \in \Par_n $ 
the following tableau set
\begin{align}{\label{blocktableaux}} 
{\cal T}_{\lambda} := 
\{  s  |  \mbox{ there is a tableau of shape } \lambda 
\mbox{ in }  [ s  ]  \}.  
\end{align} 
Let $ \left[ \lambda \right] $ be the class of $ \lambda $ under 
the equivalence relation on $ \Par_n $ 
given by $ \lambda \sim_p \mu $ if $ {\cal T}_{\lambda} =  {\cal T}_{\mu} $.
Then Murphy showed in {\it loc. cit.} that $ E_T :=\sum_{ t \in T} E_t $ and 
$ E_{[\lambda ]} := \sum_{ t \in {\cal T}_{\lambda} } E_t $ lie in $ R\Si_n$ and 
that $ \overline{E_{[\lambda ]}} \in  \F \Si_n $ 
is the block idempotent for the block given by $ \left[ \lambda \right]$. In particular 
the $ \overline {E_{[\lambda ]}} \, $'s are 
pairwise orthogonal and central in $ \F \Si_n$ with sum $1$.

\medskip
Let $ {\F \Si_{n}} \mbox{-mod} $ denote the category of finite dimensional 
${\F \Si_{n}}$-modules and 
let $$ \Ind_{n-1}^{n}:{\F \Si_{n-1}} \mbox{-mod} \rightarrow {\F \Si_{n}  \mbox{-mod}}, \,
M \mapsto {\F \Si_{n}  } \otimes_{ \F \Si_{n-1} } M $$ be the induction
functor from $ {\F  \Si_{n-1}} \mbox{-mod} $ to $ {\F  \Si_{n}} \mbox{-mod} $. 

\medskip
Assume that $ \lambda \in \Par_{n-1} $ is a subpartition of $\mu \in \Par_{n} $ and that $ \mu \setminus \lambda $ 
consists of one node of residue $ i $. Then Robinson's $i$-induction functor 
$ f_i $ is defined as 
$$ f_i :{\F \Si_{n-1}} \mbox{-mod} \rightarrow {\F \Si_{n}  \mbox{-mod}}, \, M \mapsto 
\overline{E_{[\mu ]} } \,{\F \Si_{n}  } \otimes_{ \F \Si_{n-1} } M. $$
Consider the following set $ {\cal T}_{i,n}$ of 
tableaux classes of $n$-tableaux
$$ {\cal T}_{i,n}:= \{ \, [t] \, | \,  
 s[n] = i \, \mbox{mod} \, \,  p \,\mbox{ for some (any) } s \in [t]\} $$
and set $ \overline{E_{i, n }} := \sum_{ T \in {\cal T}_{i,n} } \overline{E_T} $.
Then $ {\cal T}_{i,n}$ is a union of tableaux classes
and so 
$ \overline{E_{i,n} }$ is an idempotent in $ \F \Si_n $ and moreover $ \sum_i \overline{E_{i ,n }} = 1 $. 
We now have the following Lemma.
\begin{lemma}{\label{robinson}}
Suppose that $ M $ lies in the $[\lambda]$-block of $ \F \Si_{n-1} $. Then there is an 
isomorphism of $ \F \Si_{n} $-modules
$$ f_{i} M  \cong { \F S }_{n} \overline{E_{ i,n} }
\otimes_{\F \Si_{n-1}} M. $$
\end{lemma}
\begin{pf*}
Since the $ E_t$'s sum to $1 $, we have  
that $ {E_{[\mu ]} } $, viewed as an element of $ R \Si_n $, is 
the sum of all $ E_s $ where $ s $ is obtained from a tableau in $ {\cal T}_{\lambda} $ by adding 
an addable $i$-node. From this we deduce
$$ \overline{E_{[\mu ]}} = \overline{E_{[\mu ]}} \,  \overline{E_{[\lambda ]}} =  \overline{E_{i,n}} \ 
\overline{ E_{[\lambda ]}} .  $$ 
On the other hand $ \overline{E_{[\mu ]} } $ is central in $ \F \Si_n $ and so we get
$$ f_i M \cong {\F \Si_{n}  } \overline{E_{[\mu ]}}\otimes_{ \F \Si_{n-1} } M 
\cong {\F \Si_{n}  } \overline{E_{i,n}} \,  \overline{E_{[\lambda ]}} \otimes_{ \F \Si_{n-1} } M \cong 
{\F \Si_{n}  } \overline{E_{i,n}} \,  \otimes_{ \F \Si_{n-1} } M
$$
as claimed.
\end{pf*}

\medskip
We next introduce the notation that allows us to generalize 
the Lemma to the `divided powers'. 
Assume that $ \mu \in \Par_{n, res} $ 
and that all its ladders 
$ {\cal L}_k $ are of length $ | {\cal L}_k|  $ strictly less than $p$. 
The partition $ \mu = (6, 5, 3,1 ) $ considered above violates this condition, 
since for example $ | {\cal L}_5|  = 3 $,  
whereas the partition
$ \nu = (4,4, 3,1) $ meets it. Its $p$-residue diagram is 
$$ \nu_{res}= 
{\small{\tableau[scY]{0,1, 2, 0| 2, 0, 1 | 1, 2 , 0|0 }}} 
$$
and the ladder lengths are $ 1,1,2,2,2,1,2$, all less than $3$.
We shall need this ladder condition repeatedly and therefore introduce the following definition.
\begin{defi}
A
partition $ \mu \in \Par_{res,n} $ is called ladder restricted if all its ladders 
are of cardinality strictly less than $ p$.
\end{defi}


\medskip
Let $ \xi \in \C^{\times}$ and 
let $ {\cal H}_n(\xi) $ be the (specialized) Hecke algebra of finite type $ A $, that is 
the $ \C$--algebra on generators $ \{\,  T_i \, |\,  i= 1, 2, \ldots, n-1 \} $ 
subject to the braid relations of finite type $ A $ and the quadratic relation 
$  (T_i -\xi) (T_i +1) = 0 $.

\medskip
In [Ja], James formulated a conjecture concerning the decomposition numbers for $ q $-Schur algebras.
A special case of his conjecture is the statement that   
for $ n < p^2 $ and for $ \xi $ a primitive $p$'th root of unity,  
the decomposition numbers for $ {\cal H}_n(\xi) $ and for $ \F \Si_n$ should
coincide, via a modular reduction procedure sending $ \xi $ to $ 1 $.
We shall refer to this last statement as James' conjecture.

\begin{lemma}
Assume that  
$ n < p^2 $.  Then any $ \mu \in \Par_{res,n} $ is ladder restricted.
\end{lemma}
\begin{pf*}
Suppose that $ \mu \in \Par_{res,n}$.
Let $ {\cal L} =  {\cal L}_k$ be a ladder for $ \mu $ 
with top node $ (a_1, b_1) $ and bottom node $ (a_2, b_2) $ 
and suppose by contradiction 
that it has length $ l \ge p $ that is $ l = a_2 -a_1 +1 \ge p $.
Let $ \mu_{ \le {\cal L}} $ be the subdiagram of $ \mu $ 
consisting of the nodes $(a, b)$ satisfying $a_1 \le  a \le a_2$  and $b_2 \le b \le b_1$. It is  a
$p$-core partition.

\medskip

Since $ \mu \in \Par_{res,n}$, we know that $ \cal L $ is unbroken. 
Hence the number of nodes lying between $ (a_1, b_1) $ and $ (a_2, b_2) $ 
and on the lower border of $ \mu_{ \le {\cal L}} $ is equal to $p(p-1)+1 = p^2 - p +1$. 
On the other hand, since $ p > 2 $ there are at least $ p-1 $ more nodes in the top row of 
$ \mu_{ \le {\cal L}} $. Hence $ |  \mu_{ \le {\cal L}}  | \ge  p^2 $ and so also $ n \ge p^2 $ as required.
\end{pf*}

{\it For the rest of this section we fix $ \mu \in \Par_{res,n} $ a ladder restricted partition.}

\medskip
Let the ladders of $ \mu $ be $ {\cal L}_1, \ldots, {\cal L}_m$ and 
let  the residue of any of the nodes of $  {\cal L}_k $ be $ \iota_k $.
With respect to $ \mu $ we define numbers $ n_k $ as in ({\ref{ladder-limits}}) 
and write $ \Ind_{n_{k-1}}^{n_{k}} $ for the induction functor from the
category of finite dimensional 
$ \F \Si_{n_{k-1}} $-modules 
to the category of finite dimensional $  \F \Si_{n_{k}} 
$-modules, that is 
$$ \Ind_{n_{k-1}}^{n_{k}}: M \mapsto {\F \Si_{n_{k}}}     \otimes_{ {\F \Si}_{n_{k-1}}} M.$$
Generalizing $ {\cal T}_{i,n} $, we introduce the following set $ {\cal T}_{ { \cal L}_{k}} $ of 
tableaux classes for $\Si_{n_{k}} $  
$$  {\cal T}_{{\cal L}_{k}}
:= \{ \, [T] \, | \,  
 t[j] = \iota_k \,\,  \mbox{mod} \, \,  p \,\mbox{ for  } t \in [T]  \mbox{ and } j = 
n_{k-1} +1, \ldots, n_{k} \,\}. $$
This gives rise to the following idempotents 
$$
\begin{array}{lr} 
E^{ {\cal L}_{k} } := \sum_{ T \in {\cal T}_{ {\cal L}_{k}}} E_T  \in R  \Si _{n_{k}},  \, \, \, \, & 
\overline{E^{ {\cal L}_{k} }}  \in  \F  \Si_{n_{k}}.
\end{array} 
$$
Now since we are assuming $ | {\cal L}_k | < p $, we can define another idempotent
$$ e_{k} := \frac{1} {|{\cal L}_{k} |!} \sum_{ \sigma \in \Si_{ {\cal  L}_{k} }}  \sigma  \in  \F \Si_{n_k}.$$ 
We combine it with $\overline{E^{ {\cal L}_{k} }} $ to define 
$$ \overline{E^{({\cal L}_{k}) }} := \overline{E^{{\cal L}_k }}  e_k \in \F \Si_{n_k}, \, \, \, \,\, 
 \widetilde{e}_{\mu}:= \prod_k \overline{E^{({\cal L}_k) }}.
$$
Note that 
it is not obvious from the definition that $  \widetilde{e}_{\mu} $ is 
nonzero, although each of its factors is it. But the following Lemma
follows easily from [RH3].
\begin{lemma}{\label{combine-idempotent}}
 $ \overline{E^{({\cal L}_{k}) }} $ and $ \widetilde{e}_{\mu} $ are idempotents of $ \F \Si_{n_k} $
and $ \widetilde{e}_{\mu} =\overline{E_{ [ \mu_{lad} ] }} \prod_k e_k $.
\end{lemma}
\begin{pf*}
By Lemma 1 of [RH3] the two factors of $\overline{E^{({\cal L}_{k}) }} $ commute and so it is indeed an idempotent.
Moreover, we have that $\overline{E^{({\cal L}_{k}) }} $ 
commutes with $ { \F \Si }_{n_{k-1}} $ and so all factors of $ \widetilde{e}_{\mu} $ commute and it is 
also an idempotent. The last claim also follows from this.
\end{pf*}
We now get our divided power induction functor as 
\begin{align} 
f_{\iota}^{ (| {\cal L}_{k} |) }: \,  &
{  \F \Si }_{n_{k-1}} \mbox{-mod} \rightarrow 
{\F \Si }_{n_{k}}\mbox{-mod} \\
 &  M \mapsto { \F \Si }_{n_{k}}  \overline{E^{({\cal L}_{k}) }}
\otimes_{{\F \Si }_{n_{k-1}}}   M. 
\end{align} 
By Lemma \ref{robinson}, we have that if $ | {\cal L}_{k} | = 1 $ then $ f_{ \iota}^{ (| {\cal L}_{k} |) } = 
f_{i} $.

\medskip
With this at hand, we can now formulate the definition 
of the $ \F \Si_n $-module $ \widetilde{ A( \mu)}$, 
mentioned in the introduction of the paper. It is defined as 
\begin{equation}{\label{ladderproj}}
  \widetilde{A(\mu)}:=
 f_{\iota_m}^{ (|  \LL_m |) } \ldots f_{\iota_2}^{ (| \LL_2 |) }  f_{\iota_1}^{ (| \LL_1|) } \,   \F.  
\end{equation}
The functors $ f_{\iota_k}^{ (|  \LL_k |) } $ map projectives to projectives, and so $  \widetilde{A(\mu)} $ 
is a projective $ \F \Si_n$-module.
The following Theorem contains the basic properties of
$  \widetilde {e}_{\mu}  $ and $   \widetilde{A(\mu)} $ that shall be used throughout the paper.
\begin{Thm}{\label{expansion}}
Recall that $ \mu \in \Par_{res,n} $ is 
ladder restricted. The following statements hold: 
\newline
a) There is an isomorphism of $ \F \Si_n $-modules 
$\widetilde{A(\mu)} \cong \F \Si_n \widetilde {e}_{\mu}.$
\newline
b) Define $ {M} :=   {E_{ [ \mu_{lad} ] }  S(\lambda)} $
where $ \lambda \in \Par_n$ and $ \mu \in \Par_{res,n}$.
Then $ \overline{M} $ is a free $ \F \Si_{lad, \mu} $-module.  
In particular, if $ M \not= 0 $ then $ \widetilde{e}_{\mu} \overline{ S(\lambda)} \not=0 $.
\newline
c) For $ \lambda, \mu \in \Par_{res,n}$ there are nonnegative integers $m_{\lambda \mu}  $ 
and a triangular expansion of the form
$$ \widetilde{A(\mu)} = P(\mu) \oplus \bigoplus_{\lambda, \lambda \rhd \mu } 
\, P(\lambda)^{\oplus m_{\lambda \mu}}. $$ 
In particular $ \widetilde {e}_{\mu} \neq 0 $. 
\end{Thm}
\begin{pf*}
By definition $ \widetilde{A(\mu)} $ 
is isomorphic to 
$$ \begin{array}{l}
 \F  { \Si }_{n_{m}} \overline{E^{( {\cal L}_{m} ) }}
\otimes_{\F{ \Si  }_{n_{m-1}}} \ldots  \otimes 
\F {\Si}_{n_2} \overline{E^{( {\cal L}_{2}) }}
 \otimes_{\F{ \Si }_{n_{2}}}  
\F { \Si }_{n_{1}} \overline{E^{({\cal L}_{1}) }}
 \otimes_{ \F { \Si }_{n_{1}}}  \F.
\end{array}
$$
Note that $ n_1 = 1 $, $ \F { \Si }_{1} = \F$ and $ \overline{E^{({\cal L}_{1})}} = 1 $. 
Since $\overline{E^{({\cal L}_{k}) }} $ commutes with $ { \F \Si }_{n_{j}} $ for 
all $ j < k $, this  
simplifies to $$\F{ \Si }_n  \prod_k \overline{E^{({\cal L}_{k}) }}  = \F{ \Si }_n \widetilde {e}_{\mu}$$
proving a).
\medskip

In order to show b), we first note that $ \overline{M} $ indeed is an 
$ \F \Si_{lad, \mu} $-module, since the elements of $ \Si_{lad, \mu} $ commute 
with $ \overline{E_{  [ \mu_{lad}] }} $ 
by Lemma 1 of [RH3]. Consider now the 
set of tableaux $ T_{\mu \lambda} $ as in Lemma \ref{laddergroupaction}.
Let $ t_1, t_2 \ldots, t_k \in T_{\mu \lambda} $ 
be the maximal elements of the 
$ \Si_{lad, \mu}$ orbits in $ T_{\mu \lambda} $. 
By Lemma \ref{laddergroupaction}, 
the orbits $  \Si_{lad, \mu} \, t_i  $ are all of cardinality $ | \Si_{lad, \mu} | $. 
For each $ i $, we now check that 
the homomorphism 
\begin{equation}{\label{isinjective}} \varphi_i:
\F \Si_{lad,\mu} \rightarrow \overline{M}, \, \, \sigma \mapsto 
\overline{E_{  [ \mu_{lad}] }} \, \sigma x_{ t_i, \lambda}
\end{equation}
is injective. First of all, for $ \sigma \in \Si_{lad,\mu}$ we have 
that 
\begin{equation}{\label{modulo higher terms}}
 \overline{E_{ [ \mu_{lad} ] }} \sigma  x_{ t_i ,\lambda} =
\sigma x_{ t_i, \lambda}  = x_{ \sigma t_i, \lambda} 
\end{equation}
modulo higher terms, that is modulo an $\F$-linear combination of terms $ x_{ s \lambda} $ 
satisfying $  s  \rhd \sigma t_i $ and terms $ x_{ s t } $ satisfying that 
$ Shape(s) = Shape(t)  \rhd \lambda$.
Indeed, $ x_{ \sigma t_i, \lambda}  $ 
is an element of Murphy's standard basis 
and so the claim follows from the fact that the 
$L_i$'s act upper triangularily on the standard basis elements by Murphy's theory, 
see for example [Ma]. 
To show injectivity of $ \varphi_i$, we now suppose that 
$  \sum_{\sigma \in  \Si_{lad, \mu} }  \lambda_{\sigma} \sigma \neq 0$, and choose
$ \sigma $ with $\sigma \mu_{lad}   $ 
minimal subject to $ \lambda_{\sigma} \neq 0$. By the previous remark we find 
that the coefficient of $x_{ \sigma t_i  ,\lambda}  $ in  $ \varphi_i (\sum_{\sigma \in 
 \Si_{lad,\lambda} }  \lambda_{\sigma} \sigma )  $ is nonzero, and so 
$ \varphi_i $ indeed is injective. 

\medskip
On the other hand, the tableaux 
 $ \sigma t_i $ that appear in (\ref{modulo higher terms}), 
where $ \sigma \in \Si_{lad, \mu} $ and 
$i \in \{1, \ldots, k\} $,
are precisely those of $ T_{\mu \lambda} $, and so the elements of 
(\ref{modulo higher terms}) form a basis for $ \overline{M}$, see equation (2.4) of [Mu83]. We now 
conclude that $ \overline{M} = \oplus_i im \, \varphi_i $ and so b) is proved.

\medskip
Assume now 
that $ P(\lambda) $ is a summand of $  \widetilde{A(\mu)} $. Then 
$ \Hom_{   \F \Si_n } (   \widetilde{A(\mu)} , D(\lambda) ) \neq  0 $ and hence
$ \Hom_{   \F \Si_n } (  \widetilde{A(\mu)} , \overline{S(\lambda)} ) \neq  0 $
since $ D(\lambda) $ is a quotient of $ \overline{S(\lambda)}$ and 
$ \widetilde{A(\mu)}  $ is projective. On the other hand, by 
the definition of $ \widetilde{A(\mu)} $ we have that
\begin{equation}{\label{Hom}} 
\Hom_{   \F \Si_n } ( \widetilde{A(\mu)}, \overline{S(\lambda)} ) =  \widetilde{e}_{\mu} 
\overline{S(\lambda)}  =
\prod_k \overline{ e_k \, E_{ [ \mu_{lad} ] } \, S(\lambda)}. 
\end{equation}
We now show that $ \overline{E_{  [\mu_{lad}] } S(\lambda)} \neq 0 $ implies that $ \lambda \unrhd \mu $.
We view $ E_{ [ \mu_{lad} ]}  $ as an element of $ \Q \Si_n$ and get 
via Lemma \ref{ladder}, that in the expansion of it as a sum of $ E_t$,
only those $t$ with
$Shape(t) \trianglerighteq \mu $ can appear. On the other hand, over $ \Q $ the
standard basis $ \{ x_{s \lambda} , s \in \Std(\lambda)  \}$ for $ S(\lambda)_{\Q} := S(\lambda) \otimes_R \Q $ may be
replaced by the seminormal basis $ \{ \xi_{ s \lambda},\in
\Std(\lambda) \}$, as defined in [Mu92] via $ \xi_{ s \lambda} = E_s x_{s \lambda} $, 
and since 
$ E_t \xi_{ s \lambda} \neq 0 $ implies $ Shape(t) = \lambda $ we get the triangularity property of b). 

\medskip
To show that $ P(\mu) $ occurs
with multiplicity one in $ \widetilde{A(\mu)}  $, we set $ \lambda
= \mu $ in (\ref{Hom}) and verify that 
$  \widetilde{e}_{\mu}  \overline{S(\lambda)} $ has dimension
one over $\F$. We consider once again $ M := E_{  [ \mu_{lad}] } S(\mu) $. 
It is a free $ R $-module being a submodule of $S(\mu) $.
Let us determine its rank by extending scalars from $ R $ to 
$ \Q$. 
By Lemma \ref{ladder}, we get that in the expansion of 
$ E_{  [ \mu_{lad} ]} $ as a sum of $ E_t $'s, the occurring $t$
with $ Shape(t) = \mu $ are exactly those of the form $\sigma  \mu_{lad}  $
where $ \sigma \in \Si_{lad,\mu}$ and hence, over $ \Q $, we get a basis for $ M $ 
consisting of $ \{ \xi_{ s \mu} \} $ where $ s = \sigma \mu_{lad}   $. In other words, 
$ M $ has dimension $|  \Si_{lad,\mu} |$ over $\Q$. Hence $ M $ has rank $|  \Si_{lad,\mu} |$ over 
$ R $ and so 
$ \overline{M} $ is a free rank one $ \F \Si_{lad, \mu}$-module, by b) of the Theorem. 
Finally, we use that $ E_{  [ \mu_{lad}] } S(\mu) = E_{  [ \mu_{lad}] } D(\mu) $, 
as one gets by combining Lemma 3.35 and 3.37 of [Ma], and 
c) follows.
\end{pf*}

\medskip
\noindent
{\bf Remark.} In the paper 'Imaginary Schur-Weyl duality' [KM] that appeared on the arXiv some time after our work, 
Kleshchev and Muth study an imaginary tensor space $ M_n$ and show in their Theorem 4.2.1 that its isomorphism 
algebra is isomorphic to the group algebra of the symmetric group. It would be interesting to investigate the possible 
relationship between our work and theirs.

\section{The conjecture and the LLT-algorithm}
Let us recall the Fock space $ {\cal F}_q $ associated with the 
representation theory of the Hecke algebra $ {\cal H}_n(\xi) $ at a $ p$'th root of unity.
As a $\C(q)$-vector space, we have 
$${\cal F}_q: = \bigoplus_{\lambda \in \Par} \, {\mathbb C }(q) \lambda$$
where $ \Par := \bigcup_{n=0}^{\infty} \Par_n$ 
with the convention that
$ \Par_0 := \{ \emptyset \}$.
It is an integrable 
module for  
the quantum group $ {\cal U}_q ( \widehat{\frak sl}_p) $, 
where we use the version of $ {\cal U}_q ( \widehat{\frak sl}_p) $
that appears for example in [LLT]. This is 
the $ \C(q) $-algebra on generators 
$ e_i, f_i, \, i=0,1, \ldots, p-1 $ and $ k_h  $ for $  h $ belonging to the Cartan subalgebra $ \frak h $
of the associated Kac-Moody algebra, all 
subject to certain well known  
relations that we do not detail here. Let us explain the action of 
$ {\cal U}_q ( \widehat{\frak sl}_p) $ in $ {\cal F}_q$.
Assume that $ \gamma = \mu \setminus \lambda $ is a removable $i$-node of $ \mu $. 
We then define
\begin{align*}
N_i^l(\gamma):= & | \{ \mbox{addable } i \mbox{-nodes to the left of } \gamma\} | -
| \{ \mbox{removable } i \mbox{-nodes to the left of } \gamma\} |  \\
N_i^r(\gamma):= & | \{ \mbox{addable } i \mbox{-nodes to the right of } \gamma\} | - 
| \{ \mbox{removable } i \mbox{-nodes to the right of } \gamma\} |.
\end{align*} 
The action of $ e_i, f_i, i=0,1,\ldots, p-1$ on $ {\cal F}_q$ is now given by the following formulas
\begin{equation}{\label{action}}
f_i \lambda =  \sum_{ \mu \in \Par_n, \gamma= \mu \setminus \lambda }
q^{  N_i^l(\gamma) } \mu, \, \, \, \, \, \, \, \, 
e_i \mu =  \sum_{\lambda \in \Par_{n-1}, \gamma= \mu \setminus \lambda} q^{ - N_i^r(\gamma) } \lambda
\end{equation}
where $ \gamma $ runs over addable $ \lambda$-nodes in the first sum, and 
over removable $ \mu$-nodes in the second sum.
Note that since [LLT] use the duals of our Specht modules, 
the formulas for the action on $ {\cal F}_q$ that appear there are slightly different. 
There are similar formulas for the action of the other generators, but we leave them out.

\medskip
For $ k \in \Z $ we
let $ [k]_q:=\frac{q^k - q^{-k}}{q- q^{-1} } $ be 
the usual Gaussian integer, with the convention $[0]_q =0$, and 
define $[k]_q! := [k]_q [k-1]_q \ldots [1]_q $ 
and the divided powers 
$ f_i^{(k)} := \frac{1}{ [k]_q!} \,f_i^{k} $ and $ e_i^{ (k)} := \frac{1}{ [k]_q!} \, e_i^{k} $. 
We then introduce $ {\cal U}_{\Q} $ as the $ \Q[q, q^{-1}] $-subalgebra 
of $ {\cal U}_q (  \widehat{\frak sl}_p )$
generated by $e_i^{ (k)},f_i^{ (k)}, k=1,2,3,\ldots $. We define 
$ {\cal M}_q :=  {\cal U}_q (  \widehat{\frak sl}_p ) \, \emptyset $ and 
${ \cal M}_{\Q} :=  {\cal U}_{\Q} \,\emptyset$.

\medskip
$ {\cal M}_q $ is the basic 
module for $ {\cal U}_q (  \widehat{\frak sl}_p ) $. It is irreducible and 
therefore provided with a 
canonical basis/global crystal by Lusztig and Kashiwara's general theory.
To be more precise, let $ u \mapsto \overline{u} $ be the usual bar involution of 
${\cal U}_{\Q}  $, 
satisfying $ \overline{q} = q^{-1}, \, \overline{q^h} =q^{-h}, \,  \overline{f_i^{(k)}} = 
f_i^{(k)} $ and $ \overline{e_i^{(k)}} = 
e_i^{(k)} $ for all $ h $ and relevant $ i$. 
It induces an involution $ m \mapsto \overline{m} $ of $ {\cal M}_q$, 
satisfying $ \overline{\emptyset} = \emptyset $ and $ \overline{u v } = \overline{u}
\overline{v}$
for $ u \in {\cal U}_{\Q}  $ and 
$ v \in  {\cal M}_{\Q} $.
 
Let $ A := \{ f(q)/g(q) | f(q),g(q) \in  \Q[q], g(q) \neq 0 \} $. Then $ A $ is a local subring of 
$ \Q(q) $ with maximal ideal $ q A $ and we define $ L $ as the $A$-sublattice of $ {\cal F} $ generated by all $ \lambda \in \Par$. 
The following Theorem follows from Kashiwara and Lusztig's general theory.
\begin{Thm}
There is a unique $ \Q[q, q^{-1}] $-basis $ \{ G(\mu) \,| \, \lambda \in \bigcup_n \Par_{res,n} \}  $
for $ { \cal M}_{\Q} $, called the lower 
global crystal basis, satisfying
$$ a)\,\, G(\mu) \equiv \lambda \mbox{ mod } \, q L, \, \, \,  \, \, \, \, \, \, 
\, \, \, \, \,  \, \, \, \, \, \, \, \, \, 
b)\,\, \,  \overline{G(\mu)} = G(\mu). $$
\end{Thm}

Recall now that 
Lascoux, Leclerc and Thibon introduced in [LLT] for $ \mu \in \Par_{res,n} $ an element $ A(\mu) $ of ${ \cal M}_q$ 
called  
`the first approximation to $ G(\lambda) $'. It is defined as 
\begin{equation}{\label{firstapprox}}
  {A(\mu)}:=
 f_{\iota_m}^{ (|  \LL_m |) } \ldots f_{\iota_2}^{ (| \LL_2 |) }  f_{\iota_1}^{ (| \LL_1|) } \,   \emptyset
\end{equation}
where $ {\cal L}_1, {\cal L}_1, \ldots, {\cal L}_m $ still are the ladders for $\mu$ with 
residues $ \iota_1,  \ldots, \iota_m$. 
Based on this, 
they explain a recursive algorithm, the LLT-algorithm, that 
determines $ G(\mu)$ in terms of $ A(\lambda) $ where $ \lambda \in \Par_{res,n} $ and 
$ \lambda \unrhd \mu$. 
The following is an immediate consequence of that algorithm.
\begin{Thm}{\label{expansionGlobal}}
For $ \mu, \lambda \in \Par_{res,n} $ there is an expansion of the form 
$$ {A(\mu)} = G(\mu) + \sum_{\lambda, \lambda \rhd \mu } 
\, n_{\lambda \mu}(q) G(\lambda) $$ for certain $n_{\lambda \mu}(q) \in \Z[q,q^{-1}] $
satisfying $ \overline{n_{\lambda \mu}(q)} = n_{\lambda \mu}(q) $.
\end{Thm}
The main purpose of our paper is to study the following conjecture. 
\begin{conj}{\label{mainconjecture}}
Suppose that $ n < p^2 $ and that $  \lambda, \mu \in \Par_{res,n} $. Then 
$$ n_{\lambda \mu}(1)  = m_{\lambda \mu} $$
where $ n_{\lambda \mu}(q)$ is as in Theorem {\ref{expansionGlobal}} and
$ m_{\lambda \mu}$ as in Theorem {\ref{expansion}}.
In particular $ n_{\lambda \mu}(1) $ is a nonnegative integer.
\end{conj}

\medskip
\noindent
{\bf Remark.} 
Our main interest in studying the conjecture comes from Theorem {\ref{mainJames}}
below, which shows that it implies James's conjecture, by which we mean that 
the decomposition numbers for $ \F \Si_n $ and $ {\cal H}_n(\xi) $ coincide. 
As already mentioned in the introduction, Williamson has announced counterexamples 
to this statement within the region $n < p^2 $ suggested by James in [Ja], and so 
we refrain from speculating on the region of 
validity of the conjecture.
In the next section we do give strong experimental evidence in 
favor of the conjecture, 
see Theorem {\ref{gap}} below. Still, this evidence does not approach the 
order of the smallest counterexample given by Williamson.

We note at this point that in the cases that are covered by Theorem {\ref{gap}},  
we always have that $n_{\lambda \mu}(q) =  n_{\lambda \mu}(1)$, that is $ n_{\lambda \mu}(q)  $
is a constant polynomial, 
and so the condition $n_{\lambda \mu} =  n_{\lambda \mu}(1)$
may be necessary for the conjecture to be valid.
On the other hand, 
experimental evidence beyond the cases that are covered by Theorem {\ref{gap}} suggest 
that, as predicted by the conjecture, 
$n_{\lambda \mu}(1) $ is nonnegative even when $ n_{\lambda \mu} $ is nonconstant, 
for instance for $ p = 5 $ we have checked that $ n_{\lambda \mu}(1) \geq 0$ for all $ n < 5^2$.

\medskip
\noindent

Let $ {\cal G}(n) $ be the Grothendieck group 
of finitely generated $ \F \Si_n$-modules, and let ${\cal K}(n) $ be 
the Grothendieck group of finitely generated projective $ \F \Si_n$-modules. 
If $ M$ is a (projective) $ \F \Si_n$-module, we denote by $[M] $ its image in ${\cal G}(n) $ (${\cal K}(n) $).
We have that $ {\cal G}(n) $ and $ {\cal K}(n) $ are free Abelian groups with bases given by $ \{ [ D(\mu)] \} $ 
and $ \{ [ P(\mu)] \} $ for $ \mu \in \Par_{res,n} $. There is a non-degenerate bilinear pairing $ ( \cdot, \cdot ) $ 
between $ {\cal G}(n) $ and $ {\cal K}(n) $ which is given by $ ( [P], [M)] ) = 
\dim \Hom_{\F \Si_n} (P, M) $.
Using it, we have the following formula for the decomposition number for $ \F \Si_n$
$$ d_{\lambda \mu} = ([P(\mu)], [\overline{S(\lambda)}]).$$
These constructions and definitions can also be carried out for the Hecke algebra $ {\cal H}_n(\xi) $, 
and we shall in general use a superscript 'Hecke' for the corresponding quantities.

\medskip
Our interest in Conjecture \ref{mainconjecture} 
comes from the following Theorem.
\begin{Thm}{\label{mainJames}}
Suppose that Conjecture \ref{mainconjecture} is true. Then James' conjecture holds, 
that is $ d_{\tau \mu} = d_{\tau \mu}^{hecke} $ for all $ \tau \in \Par_n$ where $ \mu $ is as in the Conjecture.
\end{Thm}
\begin{pf*}
For any $ v \in {\cal F}_q $, 
we define $ v_{\lambda} \in \C(q) $ as the coefficient of $ \lambda  $ 
in the expansion of $ v $ in the natural basis $ \Par $ for $ {\cal F}_q$.
We first check that 
\begin{equation}{\label{weightspace}}
(A(\mu)_{\lambda})(1) = \dim_{\F} (\widetilde{e}_{\mu} \overline{S(\lambda)}).
\end{equation}
To calculate $(A(\mu)_{\lambda})(1) $ we put $q = 1 $ in the 
formula (\ref{action}) to arrive at 
$$ (f_{i_m}^{ (|  \LL_m |) } \ldots f_{i_2}^{ (| \LL_2 |) }  f_{i_1}^{ (| \LL_1|) } \, \emptyset)_{\lambda}(1). $$
But this 
is exactly the number of tableaux in $ [\mu_{lad}] $ of shape $ \lambda$, 
that is the cardinality of $ T_{\mu \lambda} $ from Lemma {\ref{laddergroupaction}}, 
as can be seen from the combinatorial description 
of $ [\mu_{lad}] $ given in that Lemma {\ref{laddergroupaction}}.
On the other hand we have 
$$ E_{[ \mu_{lad}] } S(\lambda) \, \, \,  = \sum_{ t \in [ \mu_{lad} ] } E_t S(\lambda) =  
\sum_{ t \in T_{\mu \lambda} } E_t S(\lambda). $$
But as mentioned already in the proof part b) of Theorem {\ref{expansion}}, 
Murphy gave in 
[Mu83] a basis for this space, from which we deduce that its dimension 
is the cardinality of $ T_{\mu \lambda} $, as well.
Finally we obtain ({\ref{weightspace}}), using that $\Si_{lad, \mu}$ acts faithfully 
on $ T_{\mu \lambda} $, as shown in 
Lemma {\ref{laddergroupaction}}, 
combined with part b) of Theorem 
{\ref{expansion}}.

Let us now assume that Conjecture \ref{mainconjecture} holds and let $ (a_{\lambda \mu}) :=  (n_{\lambda \mu}(1))^{-1} $.
Then we have the following formulas
\begin{equation}{\label{inversion}}
 {G(\mu)}(1) = A(\mu)(1) + \sum_{\lambda, \lambda \rhd \mu } 
\, a_{\lambda \mu} A(\lambda)(1), \, \, 
{[P(\mu)}] = [\widetilde{A(\mu)}] + \sum_{\lambda, \lambda \rhd \mu } 
\, a_{\lambda \mu} [\widetilde{A(\lambda)}] 
\end{equation}
where the last equality takes place in ${ \cal K}(n)$. We get from this last formula that 
$$ d_{\tau \mu} = 
( [\widetilde{A(\mu)} ], [S(\tau)] )  + \sum_{\lambda, \lambda \rhd \mu } 
\, a_{\lambda \mu}( [\widetilde{A(\lambda)}] ,  [S(\tau)] ) $$
which, using equation ({\ref{weightspace}}) and the definition of $ (\cdot, \cdot) $, can be rewritten as 
$$ d_{\tau \mu} = 
(A(\mu)_{\tau})(1)   + \sum_{\lambda, \lambda \rhd \mu } 
\, a_{\lambda \mu} ( A(\lambda)_{\tau})(1) = (G(\mu)_{\tau})(1) $$
where we for the last equality used the first equality of (\ref{inversion}). 
Finally, by Ariki's proof of the main Conjecture of [LLT], we know that $ (G(\mu)_{\tau})(1) = d_{\tau \mu}^{hecke} $. 
The Theorem is proved.
\end{pf*}

\noindent
{\bf Remark.} It would be interesting to study the inverse implication of the Theorem, that we believe should be true as well.

\section{Partial verification of Conjecture \ref{mainconjecture}}
In this section we give a method for verifying Conjecture \ref{mainconjecture} for 
$ n $ not too big.
It is inspired by the recent theory of KLR-algebra algebras.
Let therefore $ {\cal R}_n $ be the cyclotomic KLR-algebra (Khovanov-Lauda-Rouquier) of type $ A $
over $ \F $ and let $$ (a_{ij})_{ i, j \in \F} = \left\{ \begin{array}{cl} 2 & \mbox{if  } i = j \mod{p}  \\
-1  & \mbox{if  } i = j \pm 1 \mod{p} \\ 0   & \mbox{otherwise  } \end{array} \right. $$
be the Cartan matrix of affine type $ A_{p-1}^{(1)}$. Then $ {\cal R}_n $ 
is the $ \F$-algebra on the generators 
$$ \{  e( {\bf i}) \,| \,{\bf i} \in (\F)^n \} \cup \{ y_1, \ldots, y_{n} \}
\cup \{ \psi_1, \ldots , \psi_{n-1} \}	$$
subject to the following relations 
\begin{align}
\label{aaaa} y_1 e(\textbf{i})   = 0 \mbox{ if } i_1 =     0 \mod{p}   \\
\label{bbbb}  e(\textbf{i})   = 0 \mbox{ if } i_1 \not=     0 \mod{p}   \\
\label{kl2}e(\textbf{i})e(\textbf{j}) =\delta_{\textbf{i,j}} e(\textbf{i}) \\
\label{kl3}\sum_{\textbf{i} \in (\F)^n} e(\textbf{i}) =1 \\
\label{kl4}y_{r}e(\textbf{i}) =e(\textbf{i})y_r & \\
\label{kl5}\psi_r e(\textbf{i}) =e( \sigma_{r+1} \textbf{i}) \psi_r \\
\label{kl6}y_ry_s = y_sy_r&
&   \\
\label{kl7}\psi_ry_s = y_s\psi_r&
 \mbox{ if } s\neq r,r+1 \\
\label{kl8}\psi_r\psi_s = \psi_s\psi_r&
 \mbox{ if } |s-r|>1 \\
\label{kl9}\psi_ry_{r+1}e(\textbf{i}) =\left\{ \begin{array}{l}
(y_r\psi_r +1)e(\textbf{i}) \\
y_r\psi_r e(\textbf{i})   \\
\end{array}
\right.  &\begin{array}{l}
\mbox{if   }  i_r=i_{r+1}  \mod{p}    \\
\mbox{if   }   i_r \neq i_{r+1}  \mod{p}     \\
\end{array} \\
\label{kl10}y_{r+1}\psi_re(\textbf{i}) =\left\{ \begin{array}{l}
(\psi_ry_r +1)e(\textbf{i}) \\
\psi_ry_r e(\textbf{i})   \\
\end{array}
\right. &  \begin{array}{lcc}
\mbox{if   }  i_r=i_{r+1}  \mod{p}   \\
\mbox{if   }   i_r \neq i_{r+1} \mod{p}    \\
\end{array}
\end{align}
\begin{align}
\label{kl11}\psi_r^{2}e(\textbf{i})& =\left\{ \begin{array}{l}
0   \\
e(\textbf{i})  \\
(y_{r+1}-y_{r})e(\textbf{i}) \\
(y_{r}-y_{r+1})e(\textbf{i})  \\
\end{array}
\right.   \begin{array}{l}
\mbox{if    }  i_r=i_{r+1}  \mod{p}   \\
\mbox{if    }   i_r \neq i_{r+1} \pm 1 \mod{p}     \\
\mbox{if    }  i_{r+1}= i_r+1  \mod{p}   \\
\mbox{if    }  i_{r+1}= i_r-1 \mod{p}   \\
\end{array} \\
\label{kl12}\psi_r\psi_{r+1}\psi_re(\textbf{i})& =\left\{ \begin{array}{l}
(\psi_{r+1}\psi_r\psi_{r +1} +1)e(\textbf{i})  \\
(\psi_{r+1}\psi_r\psi_{r +1} -1)e(\textbf{i})    \\
(\psi_{r+1}\psi_r\psi_{r +1} )e(\textbf{i})    \\
\end{array}
\right.   \begin{array}{l}
\mbox{if   }  i_{r+2}=i_r=i_{r+1}-1  \mod{p}   \\
\mbox{if   } i_{r+2}=i_r=i_{r+1}+1   \mod{p}     \\
\mbox{otherwise}        \\
\end{array}
\end{align}
where $ \sigma_{r+1} = (r,r+1)$ acts on $ (\F)^n $
by permutation of the coordinates $ r, r+1$. It is an important point that $ {\cal R}_n $
is a $ \Z $-graded algebra. Indeed, the conditions
$$   \begin{array}{ccccc}
       \deg e(\textbf{i})=0,  &  &  \deg y_r=2,  &  &
 \deg  \psi_se(\textbf{i})=-a_{i_{s},i_{s+1}} \end{array}
$$
for $1 \leq r\leq n$,  $1 \leq s \leq n-1 $ and $\textbf{i} \in (\F)^{n}$
are homogeneous with respect to the relations and therefore 
define a unique $\mathbb{Z}$-grading on $ {\cal R}_n $ with degree function $ \deg$.

\medskip
As a special case of the main result of Brundan and Kleshchev's important paper [BK], we have an $ \F $-algebra 
isomorphism $f:  {\cal R}_n \cong \F \Si_n $ and so we may view $ \F \Si_n $ 
as a $ \Z $-graded algebra via $f$.


\medskip
Let us now return to the situation of the previous section.
We still set $ (a_{\lambda \mu}) :=  (n_{\lambda \mu}(1))^{-1} $ 
and consider for $ \mu \in \Par_{res,n} $ the element $ {\cal P}(\mu) \in {\cal K}(n) $ given by 
\begin{equation}{\label{conjecturalproj}}
{\cal P}(\mu) := [\widetilde{A(\mu)}] + \sum_{\lambda, \lambda \rhd \mu } 
\, a_{\lambda \mu} [\widetilde{A(\lambda)} ].
\end{equation}
In order to prove Conjecture \ref{mainconjecture} we must show for all $ \mu, \tau\in \Par_{res,n} $ that 
$ ( {\cal P}(\mu) , [D(\tau) ] ) = \delta_{ \mu \tau} $ or equivalently
\begin{equation}{\label{toproveconjecture}}
\dim_{\F} (\widetilde{e}_{\mu} D(\tau) ) + \sum_{\lambda, \lambda \rhd \mu } 
\, a_{\lambda \mu}  \dim_{\F} (\widetilde{e}_{\lambda}  D(\tau)) = \delta_{ \mu \tau}
\end{equation}
since we would then have that ${\cal P}(\mu) = [P(\mu)] $. 
The number of terms in the summation of (\ref{toproveconjecture}) is relatively small, 
so in order to verify these equations, we essentially need a way of determining the
dimension of $ \widetilde{e}_{\mu}  D(\tau) $. 

\medskip
Under Brundan and Kleshchev's isomorphism $f:  {\cal R}_n \cong  \F \Si_n  $, 
the generators 
$ \{  e( {\bf i}) \} \, \cup \{ y_1, \ldots, y_{n} \}
\cup \{ \psi_1, \ldots , \psi_{n-1} \}	$ are mapped to elements of 
$ \F \Si_n$ that we denote the same way.
For instance, we know from [BK] that $ e( {\bf i }) \in  \F \Si_n$ 
is the idempotent projector on the generalized eigenspace for the 
Jucys-Murphy elements, that is
$$ e( {\bf i } )  \F \Si_n = \{ a \in  \F \Si_n \, | \, \mbox{for all} \, \, k \,\, 
\mbox{there is} \, M \,
\mbox{such that }  (L_k - i_k)^M v = 0 \} $$
and hence we get that 
$$  \overline{E_{ [ \lambda_{lad} ] }} =  e( {\bf i }_{ \, lad, \lambda}) $$
where $ {\bf i }_{ \, lad, \lambda}  $ is the residue sequence of the ladder tableau for $ \lambda $ as 
in section 2. This is the key Lemma 4.1 of [HuMa1].
Hence via Lemma {\ref{combine-idempotent}}, 
we get that $ \widetilde{e}_{\lambda} $ is the symmetrized idempotent 
projector on a generalized weight space for the Jucys-Murphy operators. 
By orthogonality of the $ e( {\bf i})  $, we deduce 
that the dimension of $ \widetilde{e}_{\lambda}  D(\tau) $ is 
equal to the $p$-rank of $  \langle \cdot, \cdot \rangle_{\tau} $ on the restriction 
to $\widetilde{e}_{\lambda } \overline{S(\tau)}$.

\medskip
Let us now turn to the 
elements $  \psi_1, \ldots , \psi_{n-1} $ in $ \F \Si_n $. 
In [BK] they are constructed as suitable  
adjustments of certain 'intertwining elements' $ \phi_1, \ldots , \phi_{n-1}  $ 
and, as a matter of fact, in this section we shall mostly focus on these intertwining elements.
In [RH3] we found a natural realization of them, completely 
within the theory of Young's seminormal form. Indeed, we have that  $ \phi_i = \sigma_i +\frac{1}{h_L}$ 
where $ \frac{1}{h_L} = \frac{1}{L_{i-1}-L_i} $ is defined in Lemma 5 of [RH3].

\medskip
Let 
$ \{ \xi_{st} | (s,t) \in \Std(\lambda)^2, \,\lambda \in \Par_n \}$
be the seminormal basis for $ \Q \Si_n$ which was introduced in section 2.  
The action of $ \Si_n $ on it is given by the seminormal form, that is by the 
formulas of Theorem {\ref{YSFMurphy}}, 
but these formulas take place in $ S(\lambda)_{\Q} $ and therefore do not immediately 
help us in the modular setting. But note that by [Mu92] we have that 
$$ \xi_{\lambda \lambda} = x_{\lambda \lambda} \mod (\Q \Si_n)^{ > \lambda} $$ 
and so, using the seminormal form on a reduced 
expression $ d(s)= \sigma_{i_1} \ldots \sigma_{i_N} $ for $ d(s)$, we can express  
the standard basis element $ x_{s \lambda} $, when viewed as an element of $ S(\lambda) $,  
as a linear combination of the seminormal basis elements $  \xi_{s\lambda} $, but with coefficients
in $ \Q$. The next Theorem is based on this idea.

\begin{Thm}{\label{intertwiner-description}}
Let $ T $ be a tableau class and suppose that $ x \in E_T S(\lambda) $.
Then $ x $ can be written as $ x = \sum_{ t \in T} a_t \xi_{t \lambda} $.
The action of the intertwiner $ \phi_i $ is given by 
$ \phi_i x =  \sum_{ t \in T} a_t  \phi_i^{m}  \xi_{t \lambda} $ where 
\begin{equation}{\label{mysf}}
\phi_i^{m} \xi_{s \lambda}:= \left\{ 
\begin{array}{ll}
0 & \mbox{ if } | h | = 1 \, \,\,   \\
 \xi_{t \lambda}    & \mbox{ if } h > 1   \mbox{ and } s \nsim_p t  \\
\frac{ h^2-1}{h^2} \,  \xi_{t \lambda}    & \mbox{ if } h < 1   \mbox{ and } s \nsim_p t  \\
(1-\frac{1}{h}) \, \xi_{s \lambda} + \xi_{t \lambda}   & \mbox{ if } h > 1 \mbox{ and }  s \sim_p t  \\
(1-\frac{1}{h}) \, \xi_{s \lambda} + \frac{ h^2-1}{h^2} \, \xi_{t \lambda} & \mbox{ if } h < -1  \mbox{ and }  s \sim_p t   \\
\end{array}
\right.
\end{equation}
for $ t:= \sigma_i s $. We say that the first three cases of these formulas 
are the `regular' ones whereas the last two cases are the `singular' ones.
\end{Thm}
\begin{pf*}
The first statement is a consequence of the realization of the tableau class idempotent $ E_T $ 
as $ E_T = \sum_{ t \in T} E_t $ and the fact that $ E_t S(\lambda)_{\Q} = \Q  \xi_{t \lambda} $
for $ t \in \Std(\lambda) $.

\medskip
In order to prove the second statement, we need to recall the construction 
of $ \phi_i$ from [RH3]. Let $ S:= [s], T:=[t]   $ be as in the announcement of 
the Theorem and suppose first that $ S \neq T $. 
Choose arbitrarily $ t \in T $ and define $ c_T(i-1)  := c_t(i-1) \in R$ and  
$ c_T(i)  := c_t(i) \in R$ and set $ h_T(i):= c_T(i-1) - c_T(i) $. 
Although $ h_T(i) \in R $ depends on the choice of $ t \in T$, we showed in [RH3] that 
for any 
$ a \in  { E_T S(\lambda)} $, we have 
that $ (L_{i-1} - L_i - h_T(i))^N a$ belongs to $ p{ E_T S(\lambda)} $ for $ N $ sufficiently big, 
independently of the choice of $ t$. 
Then $ \frac{ 1} { L_{i-1} - L_i} $ is the linear transformation on 
$\overline{ E_T S(\lambda)} $
given by the corresponding geometric series. To be precise, for 
$ a \in \overline{ E_T S(\lambda)} $ it is given by  
\begin{equation}{\label{powerseries}}
 \frac{ 1} { L_{i-1} - L_i} \, a := \frac{1}{h_T(i)} 
\sum_k (-1)^{k}  \left(\frac{L_{i-1} - L_i -  h_T(i) }{h_T(i)} \right)^k a
\end{equation}
where the sum may be assumed to be finite by the above remark.
Finally, $ \phi_i $ is the linear transformation on 
$\overline{ E_T S(\lambda)} $ given by
$ \phi_i :=\sigma_i + \frac{ 1} { L_{i-1} - L_i} $. Note the slight 
variation from [RH3], where we used the definition
$ \phi_i :=\sigma_i - \frac{ 1} { L_{i-1} - L_i} $. With this convention, 
$ \phi_i $ coincides exactly with Brundan and Kleshchev's element $ \phi_i$ and 
it verifies the following intertwining property
\begin{equation}{\label{intertwinerprop}} 
E_S \phi_i  =  \phi_i E_T
\end{equation}
in 
$  \Hom_{\F} (E_T (\F \Si_n),  \F \Si_n) $, corresponding to Lemma 7 of [RH3].

\medskip
Now, recall that the argument in [RH3] to show that
$ (L_{i-1} - L_i - h_T(i))^N a \in p S(\lambda) $ for $ N >> 0$ used the expansion of
$ a \in E_T S(\lambda) $ in the seminormal basis $ \xi_{t \lambda} $. For each term 
of this expansion we indeed get 
$$ (L_{i-1} - L_i - h_T(i))^N \xi_{t \lambda} = 
(c_t(i-1) - c_t(i) - h_T(i))^N \xi_{t \lambda} \in p S(\lambda)  $$
for $ N >> 0$. 
We conclude from it that the series (\ref{powerseries}) can be calculated by lifting 
$ a \in \overline{E_T S(\lambda)}$ to $ a \in E_T S(\lambda)$, then expanding this $ a $ in
the $ \xi_{t \lambda} $'s, next 
applying the series to each term and finally reducing modulo $p$. On each of these terms we get 
\begin{equation}
 \frac{ 1} { L_{i-1} - L_i} \, \xi_{s \lambda} = 
\frac{1}{h_T(i)} 
\sum_k (-1)^{k}  \left(\frac{ c_{i-1}(t) - c_{i}(t) -  h_T(i) }{h_T(i)} \right)^k \xi_{s \lambda} 
\end{equation}
that equals $ \frac{ 1} { c_{i-1}(t) - c_i(t)} \, \xi_{s \lambda}$.
From this the regular cases (\ref{mysf}) of the Theorem follow by applying the classical 
formulas for Young's seminormal form, that is Theorem 
\ref{Murphyysf}. 

The singular cases are easier to handle, since 
we then have $ \phi_i = \sigma_i +1 $ and so we finish by applying Theorem 
\ref{Murphyysf} once again.
\end{pf*}

\noindent
{\bf Remark}. In their recent work [HuMa2], which is independent of ours, Hu and Mathas prove
in a systematic way results that are similar to the Theorem. They use them to 
develop the theory of KLR-algebras, including integral versions of them, 
completely from the point of view of Young's seminormal
form. 

\begin{Thm}{\label{norm}}
Let $ \lambda \in \Par_n $ and 
let $ \{ \xi_{ s \lambda} \, | \, s \in \Std(\lambda) \} $ be 
the seminormal basis for $ S(\lambda)_{ \Q} $. Then we have that 
$$ \langle \xi_{ s \lambda} , \xi_{ s \lambda} \rangle_{\lambda} = \gamma_s. $$
\end{Thm}
\begin{pf*}
This is contained in Murphy's papers where 
it is shown by induction. The induction basis is given by $ \xi_{\lambda \lambda } = x_{\lambda \lambda } $ 
and the induction step by Young's seminormal form (\ref{ysf}). 
\end{pf*}

With the above Theorems at our disposal we can now describe an algorithm for calculating 
$\dim_{\F} \widetilde{e}_{\lambda}  D(\tau) $, or equivalently  
the rank of $  \langle \cdot, \cdot \rangle_{\tau} $ on  $\widetilde{e}_{\lambda } S(\tau)$. 
Note that Step 3 of the algorithm depends on the results from our previous work [RH3], 
giving a cellular basis in terms of the $ \phi_i$'s.

\medskip
\noindent
{\bf Algorithm.} 
\newline
Step 1. Determine the set $ T_{\lambda \tau}= \{ s \in [ \lambda_{lad}] \, | \, Shape(s) = \tau \} $. 
As indicated above, $ T_{\lambda \tau}$ can be read off from
the calculation of the first approximation of $ A(\mu) $ at $ q = 1 $, since 
the successive actions of $ f_i$ may be viewed as producing tableaux rather than partitions.
\newline
Step 2. 
Write the elements of $  \{ d(s) \, | \, s \in T_{\lambda \tau}\} \subset \Si_n $ 
as reduced products of simple transpositions $ \sigma_i$.
The longest element of $ \Si_n $ has length $ l(w_0) = n(n-1)/2$, and so 
each of the reduced products has less than $n(n-1)/2$ terms.
\newline
Step 3. For each $ d(s) = \sigma_{i_k} \ldots \sigma_{i_1} $ from step 3, calculate 
$  \phi_{ i_k } \ldots \phi_{ i_1 } \xi_{\lambda \lambda } $ using 
({\ref{mysf}}). 
By Theorem 2 and Lemma 9 of [RH3], we get in this way an $R$-basis for $ e_{ [\lambda_{lad}] } S(\tau) $. 
The basis elements are given as linear combinations of seminormal basis elements.
The number of terms $ \xi_{ t \lambda} $ in this expansion will be less than 
$ 2^B $ where $ B $ is the number of indices in the reduced expression for $ d(s) $ that involve
the singular cases of ({\ref{mysf}}). 
\newline
Step 4.
Symmetrize each basis element from the previous step with respect to the ladder group 
$\Si_{lad, \lambda} $, to get a basis for $\tilde{e}_{\lambda } S(\tau)$.
\newline
Step 5.
Calculate the matrix of the form 
$  \langle \cdot, \cdot \rangle_{\tau} $ on $\tilde{e}_{\lambda } S(\tau)$ with respect to the basis 
given in the previous step. Since the basis elements are expanded
in terms of the seminormal basis, this step now follows easily from the previous Theorem.
The matrix will have values in $ R $ although the coefficients of the expansions are rational.
\newline
Step 6.
Reduce the matrix modulo $ p $ and determine its rank.

\medskip
\noindent
{\bf Remark.} 
The algorithm can also be implemented using the classical 
Young's seminormal form, that is Theorem \ref{Murphyysf}. 
On the other hand, that algorithm will be much less efficient with expansions that grow too fast. In fact, 
the main point of our algorithm, as presented above, is that the indices of 
the reduced expressions will mostly correspond to the regular cases of 
(\ref{mysf}), thus reducing, as much as possible, the doubling up of terms. 

\medskip
\noindent
{\bf Example.} 
Suppose that $ p = 3 $. We verify Conjecture \ref{mainconjecture} for 
$ \Si_5 $ using our algorithm. Although the LLT-algorithm does not involve any subtractions
in this example, the example is still big enough to illustrate our algorithm.

We have that $$ \Par_{res,5} = \{ 
\{    [3,2],   [3,1^2],   [2^2,1],   [2,1^3], 
  [1^5] \}. $$
The ladder groups are $$ \Si_{lad,   [3,2] } =   \Si_{lad,  [3,1^2] } = \{ (3,4) \}, \,
\Si_{lad,  [2^1,1] } =   \Si_{lad, [2,1^3] } =  \Si_{lad, [1^5] } =  1 . $$  
The first approximations are 
$$
\begin{array}{ll} 
A([3,2]) = [3,2] + q [4,1], & 
A( [3,1^2] ) = [3,1^2],  \\
A([2^2,1] ) = [2^2,1] + q [5], & 
A([2,1^3] ) = [2,1^3] + q [2^2,1], \\
A([1^5] ) = [1^5] + q [3,2]. 
\end{array} 
$$
From this we conclude, as already mentioned above, that $ G(\lambda) = A(\lambda) $ 
for all $ \lambda \in \Par_{res,5}$.
Thus, Conjecture \ref{mainconjecture} is in this case the affirmation that 
$ P(\lambda) = \widetilde{A(\lambda)} $ for all $ \lambda \in \Par_{res,5}$,
or by the above that 
$$ \dim \widetilde{e}_{ \lambda} D(\tau) = \delta_{ \lambda \tau} \, \, \,\, \, \mbox{  for all  } 
\lambda, \tau \in \Par_{res,5}. $$
We calculate the rank of $ \langle \cdot, \cdot \rangle_{\tau} $ on each symmetrized 
weight space $  \widetilde{e}_{ \lambda} S(\tau)$, using our algorithm.
Recall that $ \dim  \widetilde{e}_{ \lambda} S(\tau)$ can be read off from the first approximation 
$ A(\lambda) $. For example we have that 
$ \dim \widetilde{e}_{ [3,2]} S([4,1]) =1 $ since the coefficient $ q $ evaluates to
$1$, although the eigenspace $ \widetilde{e}_{ [3,2]} S([4,1])$ is irrelevant to us since
$ [4,1] \notin \Par_{res,5}$. 

\medskip
By going through the first approximations, we see that the 
relevant eigen--spaces are $ \widetilde{e}_{ [2,1^3]} S([2^2,1 ]) $ and
$ \widetilde{e}_{ [ 1^5 ]} S([3,2]) $, both of dimension one. We verify that 
the ranks of the corresponding forms are zero, or equivalently that the forms
are zero.

\medskip
We first consider $ \widetilde{e}_{ [2,1^3]} S([2^2,1 ]) $. 
The residue diagram $ res_{[2,1^3]} $ of $ \lambda := [2,1^3] $ is  
$$
\begin{array}{lll}
res_{[2,1^3]} = {\small{\tableau[scY]{0,1| 2 | 1 | 0  }}} \,\, \, \, \, \,\,  &
t:= {\small{\tableau[scY]{1, 2| 3, 5 | 4 }}} \,\, \, \, \, \,\,& 
res_{\tau } = {\small{\tableau[scY]{0,1| 2 ,0 | 1   }}} \,\, \, \, \, \,\,  
\end{array}
$$
and so we have $ {\bf{i}}_{lad, [2, 1^3]  } = ( 0,1,2,1,0)$. 
The only tableau of shape $ \tau :=[2^2, 1] $ in the ladder class of $ \lambda  $ 
is therefore $ t $ as given above. We have $ d(t) = (4,5) = \sigma_5$ and so 
we get from formula (\ref{mysf}) that the basis for $ \widetilde{e}_{ [2,1^3]} S([2^2,1 ]) $
is $ \{ \phi_5 \xi_{ \lambda \lambda } \}  =\{  \xi_{ t \lambda } \}$.
Finally, by Theorem {\ref{norm} we get that $ \langle \xi_{ t \lambda }, \xi_{ t \lambda } \rangle_{\lambda} 
= 3 = 0 \mod 3 $, as claimed.

\medskip
We next consider $ \widetilde{e}_{ [1^5]} S([3,2 ]) $ where we basically proceed as before.
The residue diagram $ res_{ [1^5]} $ is 
$$
\begin{array}{lll}
res_{[1^5]}:= {\small{\tableau[scY]{0| 2 | 1 | 0 | 2 }}} \,\, \, \, \, \,\,  &
s:= {\small{\tableau[scY]{1, 3, 5| 2, 4 }}} \,\, \, \, \, \,\,  &
res_s:= {\small{\tableau[scY]{0,1,2| 2, 0 }}} 
\end{array}
$$
and we have $ {\bf{i}}_{lad, [2, 1^3]  } = ( 0,2,1,0,2)$. 
The only tableau of shape $ \nu := [3, 2] $ in the ladder class of $ [1^5]  $ 
is $ s $ as given above. We have $ d(s) = (2, 3)  (4,5) (3, 4) = \sigma_3  \sigma_5  \sigma_4 $ and so 
we get from formula (\ref{mysf}) that the basis for $ \widetilde{e}_{ [1^5]} S([ 3,2 ]) $
is 
$$ \{ \phi_3  \phi_5  \phi_4 \xi_{ \lambda \lambda } \}  =\{  \xi_{ s \lambda } \} $$
and then by Theorem {\ref{norm} we get $ \langle \xi_{ s \lambda }, \xi_{ s \lambda } \rangle_{\lambda} 
= 3 = 0 \mod 3 $, as claimed. This concludes the verification of Conjecture \ref{mainconjecture}, 
and then by Theorem {\ref{mainJames}} also of James' conjecture, in this case.

\medskip
We have implemented the algorithm using the GAP-system and have found the following results, 
that without doubt can be improved on.  
\begin{Thm}{\label{gap}}
If $ p = 3 $ then Conjecture \ref{mainconjecture} is true. \newline
If $ p = 5 $ then Conjecture \ref{mainconjecture} is true for $ n < 16$. \newline
If $ p = 7 $ then Conjecture \ref{mainconjecture} is true for $ n < 19$. \newline
If $ p =  11 $ then Conjecture \ref{mainconjecture} is true for $ n < 22$. \newline
If $ p = 13$ then Conjecture \ref{mainconjecture} is true for $ n < 22$. \newline
\end{Thm}
\noindent
{\bf Remark.} 
The Theorem provides, via Theorem {\ref{mainJames}}, decomposition numbers 
for $ \F \Si_n$. But as pointed out to us by A. Mathas, the partitions involved
in the Theorem
are all of $p$-weight less than three and hence the corresponding decomposition numbers
are already present in the literature, see M. Richards' paper [Ri] for the weight one and two cases, and 
M. Fayers's paper [F] for the weight three case.

\section{$ \widetilde{ A( \lambda)}$ as a graded module}
In this section we show that $ \widetilde{ A( \lambda )}$ may be viewed as a graded module
for $ \F \Si_n = {\cal R}_n $. We do so by showing that $ \widetilde{e}_{\lambda} $ is 
a homogeneous idempotent in $\F \Si_n $, necessarily of degree zero.

\medskip
Let $ \widetilde{ {\cal R}}_n $ be the noncyclotomic KLR-algebra of type $ A $ 
over $ \F$, or more precisely the $ \F$-algebra 
on the same generators and relations as ${ {\cal R}}_n  $, but without 
the relations (\ref{aaaa}) and (\ref{bbbb}). This is the algebra that was first considered in [R] and 
in [KL] from a diagrammatic point of view. Assume that $ \lambda \in \Par_{res,n} $ is ladder restricted, 
with associated ladder residue sequence $ {\bf i }_{ \, lad, \lambda}  $
and let $ n_0, \ldots, n_m $ be as in ({\ref{ladder-limits}}).
Then for all $ k $, the residues $ i_{ n_{k-1} +1 }, \ldots, i_{n_k} $ are equal, being 
the residue $ \iota_k$ of the $k$'th ladder $ {\cal L}_k $ of $ \lambda_{lad} $.
Let $ \widetilde{ {\cal R}}_{k,n} $ be the subalgebra of $ \widetilde{ {\cal R}}_n $
defined by 
$$ \widetilde{ {\cal R}}_{k,n}  := \langle \, e( {\bf i }_{ \, lad, \lambda}), \,
y_i e( {\bf i }_{ \, lad, \lambda}), \psi_j e( {\bf i }_{ \, lad, \lambda}) \,
| \, a_k \le i \le b_k, \,  \, a_k \le j \le b_k-1  \,  \rangle$$
where we write $ {\cal L}_k = \{a_k, a_k+1, \ldots, b_k \} $ or just 
$ \{a, a+1, \ldots, b \} $ 
for simplicity.
Then $ \widetilde{ {\cal R}}_{k,n} $ is isomorphic to the nilHecke algebra, or
to be more precise, setting $ y_i := y_i e( {\bf i }_{ \, lad, \lambda}) $ and
$ \partial_j :=  \psi_j e( {\bf i }_{ \, lad, \lambda}) $,
to the infinite dimensional $\F$-algebra 
generated by $ y_r $ and $ \partial_s $ subject to the relations 
$$ 
\begin{array}{ll}
\partial_r^2 =0 &   \\
\partial_r \partial_{r+1}  \partial_r  = \partial_{r+1} \partial_r \partial_{r+1} &
\\ 
\partial_r \partial_{s}    = \partial_s \partial_r   &
\mbox{if } | r- s| > 1 \\ 
y_r y_s = y_s y_r & \\
\partial_r y_{r+1} - y_r \partial_r =1, & y_{r+1} \partial_r  - \partial_r y_r  =1 \\
\partial_r y_{s} =  y_s \partial_r  & \mbox{if }  s \neq r, r+1  \\ 
\end{array}
$$
with one-element $ e( {\bf i }_{ \, lad, \lambda})$. 
Note that to get the nilHecke algebra presentation used in [KL], one should use the 
isomorphism given by $  \partial_i \mapsto -\partial_i, y_i  \mapsto y_i$.
For $ w \in    \Si_{ {\cal L}_{k} } $ we define $ \partial_w = 
\partial_{i_1} \ldots \partial_{i_k}  $ where 
$ w = \sigma_{i_1}  \ldots \sigma_{i_k} $ is a reduced expression; this is 
independent of the chosen reduced expression.
Note that $  \prod_k \widetilde{ {\cal R}}_{k,n} $, 
is a subalgebra of $ \widetilde{ {\cal R}}_{n} $.
For  
$ w_{0,k} \in  \Si_{ {\cal L}_{k} } $ the longest element we define 
$$
\begin{array}{rl}
e_{KLR,\lambda, k,n}: &= (-1)^{ m_k( m_k -1) /2 } \, 
\partial_{w_{0,k}} y_{a_k}^{m_k-1} y_{a_k+1}^{m_k-2} \ldots y_{b_k-1} \in \widetilde{ {\cal R}}_{k,n} \\
e_{KLR,\lambda,n}: &=  \prod_k e_{KLR,\lambda, k,n}  \in  \widetilde{ {\cal R}}_{n}
\end{array}
$$
where $ m_k = | {\cal L}_k | $.
Then, by section 2.2 of [KL], $ e_{KLR,\lambda, k,n} $ is a 
homogeneous idempotent of $ \widetilde{ {\cal R}}_{k,n} $ and
hence $ e_{KLR,\lambda,n} $ is a homogeneous idempotent of $ \widetilde{ {\cal R}}_{n} $.
Let $ g : \widetilde{ {\cal R}}_{n} \rightarrow { {\cal R}}_{n} $ be the  
quotient map.

\medskip
In the following Theorem we show that our idempotent $  \widetilde{e}_{\lambda} $ 
is equal to $ g(e_{KLR,\lambda,n}) $. The main difficulty in showing this is prove that 
$ g(e_{KLR,\lambda,n}) \neq 0$,  since in general the $ y_i $'s are zero divisors in $  {\cal R}_{n} $.
Our solution to the problem uses once again Young's seminormal form and goes back to ideas of Hu and Mathas. In fact, 
the arguments that lead to equation ({\ref{onlyoneterm}}) below are closely related to the arguments leading to the crucial Theorem 4.14 
of [HuMa1].

\begin{Thm}{\label{idempotents}}
Assume that $ \lambda \in \Par_{res,n} $ is ladder restricted.
Then $ g(e_{KLR,\lambda,n}) = \widetilde{e}_{\lambda}$.
In particular, $ \widetilde{ A( \mu)} $ may be considered a graded module for $ \F \Si_n $.
\end{Thm}
\begin{pf*}
Clearly, $g(e_{KLR,\lambda,n})$ is an idempotent, although possibly zero. 
By the definitions we have 
$ g(e_{KLR,\lambda,n}) = \prod_k g( e_{KLR,\lambda, k,n}) $ where 
$$ g( e_{KLR,\lambda, k,n}) = (-1)^{ m_k( m_k -1) /2 } \, \partial_{w_{0,k}} y_{a_k}^{m_k-1} y_{a_k+1}^{m_k-2} \ldots y_{b_k-1}  $$
where we use the same notation for $ x \in \widetilde{ {\cal R}}_{n}$ and its image 
$ g(x) \in  {\cal R}_{n}$.  

\medskip 
Let us now recall Brundan and Kleshchev's construction of the element $ g( \partial_j )$, 
that is $ \psi_j e(  {\bf i }_{ \, lad, \lambda}) $.
As already mentioned it is an adjustment of the intertwining element $ \phi_j$, that we described in 
Theorem {\ref{intertwiner-description}}.
On the other hand, since $ a_k \le j \le b -1 $ we are in the singular case of 
Theorem {\ref{intertwiner-description}} and so we have $ \phi_j = \sigma_j +1 $
when acting in the generalized eigenspace corresponding to $ e( {\bf i }_{ \, lad, \lambda}) $.
In general, the adjustment element $ q_j( {\bf i }_{ \, lad, \lambda} ) $ 
satisfying $ \phi_j = \psi_j q_j( {\bf i }_{ \, lad, \lambda} )^{-1}$ 
is an invertible power series in the $ y_i$'s. 
Let now $ j $ be any number such that $ a_k \le j \le b_k -1 $. 
Then there is a reduced expression for $ w_{0,k}$
of the form $w_{0, k } = \sigma_j w^{\prime} $ for some $ w^{\prime}$
from which we deduce that $$ \partial_{w_{0,k}} = (\sigma_j+1) q_j( {\bf i }_{ \, lad, \lambda} )^{-1} a $$
for some $ a $. Hence we conclude that $g( e_{KLR,\lambda, k,n})$
is invariant under the action of $ \sigma_j $ for all $ j $ such that $ a_k \le j \le b_k -1 $
and so we have that $ g(e_{KLR,\lambda,n}) = c  \widetilde{e}_{\lambda}$ for some scalar $ c$.
Since $ g(e_{KLR,\lambda,n})$ is an idempotent, the only possibilities for $ c $ are now $ c = 1 $ or 
$ c=0 $.

\medskip
It remains to show that $ g(e_{KLR,\lambda,n}) \neq 0 $, or that $ c = 1$.
Let $ x_{ \lambda_{lad}, \lambda} \in S(\lambda)$ be the Murphy basis element as above. 
From equation ({\ref{isinjective}}) we have that 
$ \widetilde{e}_{\lambda} x_{ \lambda_{lad}, \lambda} \neq 0$, so in order to show that 
$ c = 1 $ it is enough to check the equality
\begin{equation}{\label{equality}} 
g(e_{KLR,\lambda,n}) \widetilde{e}_{\lambda} x_{ \lambda_{lad}, \lambda} =  
\widetilde{e}_{\lambda} x_{ \lambda_{lad}, \lambda}.    
\end{equation}
Without loss of generality, it is enough to consider the $k$'th ladder $ {\cal L}_k = {\cal L}$
and to assume that $\Si_{ lad, \lambda}=  \Si_{\cal L} $. 
Let $ w_{k,0} = w_{0} $ be the longest element of $ \Si_{\cal L} $. 
In order to show ({\ref{equality}}), we first lift both sides to $ \Q $, 
then expand 
in terms of the seminormal basis $ \xi_{ \sigma \lambda, \lambda}, \,\sigma \in \Si_{ {\cal L} } $ and finally 
verify that in both expansions the coefficient of $ \xi_{ w_{0}\lambda. \lambda} $ is the same, up to 
a unit in $ R^{\times}$. 

\medskip
We first consider the right hand side of 
({\ref{equality}}). It is the reduction modulo $ p $ of the following element, 
that exists over $ R$, but that shall also be considered over $ \Q$
$$ LIFT:=
 1/ |  \Si_{\cal L}  | \sum_{ \sigma \in  \Si_{\cal L} }
\sigma E_{ [ \lambda_{lad}] } x_{ \lambda_{lad}, \lambda}. $$
Using 
Young's seminormal form ({\ref{YSFMurphy}}) repeatedly 
we get that the coefficient of $ \xi_{ w_{0} \lambda, \lambda} $ in 
$ LIFT$ is $1/ |  \Si_{\cal L} | $.

\medskip
We next work out the left hand side of ({\ref{equality}}), using 
the realization in [BK] of the $ y_i$'s as Jucys-Murphy elements. Indeed,  
we have from (3.21) of {\it loc. cit.} that 
\begin{equation}{\label{y-realization}} 
y_r = \sum_{ {\bf i} \in ({\F})^n } (L_r - i_r) e( { \bf i } ). 
\end{equation}
For $ T $ a tableau class with residue sequence $ {\bf i} $, 
this can be lifted to $ R$ to 
\begin{equation}{\label{y-realization2}}  
y_r = \sum_{ t \in T } ( c_{t}(r) - \hat{i}_r)  E_t
\end{equation}
where $ \hat{i}_r \in \Z $ is chosen such that $ \hat{i}_r \mbox{ mod } p = i_r$. 
Writing 
$ LIFT=
\sum_{ \sigma \in \Si_{\cal L}} a_{\sigma} \xi_{\sigma \lambda_{lad} , \lambda} $ for $ a_{\sigma} \in \Q$, we get
\begin{equation}{\label{y-realization3}}  
y_{a}^{m-1} y_{a+1}^{m-2} \ldots y_{b-1} \, LIFT=
y_{a}^{m-1} y_{a+1}^{m-2} \ldots y_{b-1} \, 
\sum_{ \sigma \in  \Si_{ {\cal L}}} a_{\sigma} \xi_{\sigma \lambda_{lad} , \lambda}.
\end{equation}
Now, for the factors of $y_a^{m-1} $ 
we choose for the lift $ \hat{i}_r $ as in ({\ref{y-realization}})  
all the contents that appear in $ { \cal L} $ each exactly once, except $ c_{\lambda_{lad}}(a)$.
Similarly, 
for the factors of $y_{a+1}^{m-2} $ 
we choose for $ \hat{i}_r $ as in ({\ref{y-realization}})  
all the contents that appear in $ { \cal L}$ each once, except this time $ c_{\lambda_{lad}}(a)$ 
and $ c_{\lambda_{lad}}(a+1)$ and so on.
With these choices, only one term survives in the action of the $ y_i$'s
in ({\ref{y-realization3}}), giving 
\begin{equation}{\label{onlyoneterm}}  
y_{a}^{m-1} y_{a+1}^{m-2} \ldots y_{b-1} \, 
\sum_{ \sigma \in  \Si_{{\cal L}}}  a_{\sigma} \xi_{\sigma \lambda_{lad} , \lambda}  
=  u a_{1 } p^{m(m-1)/2}  \xi_{\lambda_{lad} , \lambda}  
\end{equation}
for some $ u \in R^{\times} $ where we use
that $ | {\cal L} | < p $ to get the factor $p^{m(m-1)/2}$. 

\medskip
We now have to calculate $ a_1 $, the coefficient of 
$ \xi_{\lambda_{lad} , \lambda}  $ in $LIFT$.
The coefficient of $ \xi_{w_0 \lambda_{lad}  , \lambda} $ is $ 1/|  \Si_{\cal L}  | $, as already mentioned above, so 
for this we need a formula relating the coefficients of the $ \xi_{s \lambda_{lad}, \lambda} $'s.
But using Young's seminormal form, Theorem \ref{Murphyysf}, on the equality
$$ \sigma_i \sum_{ \sigma \in { \Si_{{\cal L}}} } a_{\sigma} \xi_{\sigma \lambda_{lad} , \lambda} =
\sum_{ \sigma \in \Si_{ {\cal L} }} a_{\sigma} \xi_{\sigma \lambda_{lad} , \lambda} \,\,\,\,\, \mbox{    for all }i$$
we get that 
$$ a_{ \sigma  } =( (h-1)/h) \, a_{ \sigma_i \sigma  } \mbox{     if }  \sigma_i  \sigma > \sigma $$
where $ h $ is the radial distance between $ \sigma  \lambda_{lad}  $ 
and $\sigma_i \sigma  \lambda_{lad}   $ and then we 
find that $ a_1 =  p^{-m(m-1)/2} /|  \Si_{\cal L}| $.

From this we deduce, via ({\ref{onlyoneterm}}), that 
\begin{equation}{\label{y-realization3}}  
y_{a}^{m-1} y_{a+1}^{m-2} \ldots y_{b-1} \, LIFT=  u_{1 }  \xi_{\lambda_{lad} , \lambda}  
\end{equation}
where $ u_1 \in R^{\times}$.
Finally, using  
$ \partial_{j} = (\sigma_j+1) q_j( {\bf i }_{ \, lad, \lambda} )^{-1}  $ once again, 
together with Young's seminormal form, we find that 
the coefficient of $ \xi_{w_0 \lambda_{lad}  , \lambda} $ in 
$$ \partial_{w_{0,k}} y_{a}^{m-1} y_{a+1}^{m-2} \ldots y_{b-1} \, LIFT $$ 
is a unit of $ R$.
The Theorem is proved.
\end{pf*}

Via the ismorphism 
$ {\F \Si_{n}} \cong {\cal R}_n$, 
we may introduce $ {\F \Si_{n}} \mbox{-grmod} $, the category of finite dimensional 
graded $ {\F \Si_{n}} $-modules. For $ M $ an object of $ {\F \Si_{n}} \mbox{-grmod}$, 
we have a decomposition of $ \F $-spaces
$$ M = \oplus_{i \in \Z} M_i $$
such that 
$ ({\F \Si_{n}})_i M_j \subset M_{i+j}$ for all $i,j$.

\medskip
Let $ v:{\F \Si_{n}} \mbox{-grmod}  \rightarrow {\F \Si_{n}} \mbox{-mod} $ be 
the forgetful functor. Using the results from [CF] on the representation theory of general $\Z$-graded 
rings, we have that an object $ M $ in $ {\F \Si_{n}} \mbox{-grmod} $ 
is indecomposable if and only if $ v(M) $ is indecomposable in $  {\F \Si_{n}} \mbox{-mod} $
and that every indecomposable object in $  {\F \Si_{n}} \mbox{-mod} $ is of the form 
$ v(M) $ for some $ M $.
Moreover, by Lemma 2.5.3 of [BGS] we know that the grading on indecomposable modules 
is unique up to a shift in degree. That is, 
if $ M $ and $  N $ are indecomposable objects in 
$ {\F \Si_{n}} \mbox{-grmod} $ satisfying $ v(M) \cong v(N) $, then $ M \cong N\langle k \rangle $ 
for some $ k \in \Z$ where $ N\langle k \rangle $ is the degree shift of $ N $ of order $ k$, that is
$$  N\langle k \rangle_i := N_{i-k}.$$
Let us write $ \widetilde{A(\mu)^{gr}} $ for $ \widetilde{A(\mu)} $ considered as a
graded module, via the previous Theorem.
By applying the above mentioned results from the literature, we now get, 
corresponding to b) of the Theorem {\ref{expansion}}, 
a triangular expansion
$$ \widetilde{A(\mu)^{gr}} = P(\mu)^{gr}  \oplus \bigoplus_{\lambda, \lambda \rhd \mu } 
\, (P(\lambda)^{gr}\langle k_{\lambda \mu j} \rangle)^{\oplus m_{\lambda \mu}} $$
for certain integers $m_{\lambda \mu}, k_{\lambda \mu j}  $ where 
$ P(\mu)^{gr} $, satisfying $v(P(\mu)^{gr} ) = P(\mu)$, 
is chosen as the first term in the expansion 
for all $ \mu $.

\medskip
Formally, the relationship between $ \widetilde{A(\mu)^{gr}}$ and $ P(\mu)^{gr} $ 
is now the same as the one between the Bott-Samelson bimodule and the indecomposable 
Soergel bimodule in Soergel's theory of bimodules over the coinvariant ring 
of a Coxeter group, see [So]. It would therefore be interesting to investigate to what extent 
the methods of Elias and Williamson's paper [EW] for 
solving Soergel's conjecture can be applied to our situation.

\sc Instituto de Matem\'atica y F\'isica, Universidad de Talca, Chile, 
steen@inst-mat.utalca.cl.


\begin{thebibliography}{99}
\bibitem[A]{Ariki}
S. Ariki, {\it On the decomposition numbers of the Hecke algebra of $G(m, 1, n)$}, J.
Math. Kyoto Univ. {\bf 36} (1996), 789--808.
\bibitem[BGS]. A. Beilinson, V. Ginzburg, W. Soergel, {\it Koszul duality patterns in representation theory}, 
J. Amer. Math. Soc. {\bf 9} (1996), 473--527.
\bibitem[BF]{brundan-kleshchev}
J. Brundan, A. Kleshchev, {\it Blocks of cyclotomic Hecke algebras and Khova-nov-Lauda algebras}, 
Invent. Math. {\bf 178}  (2009),  no. 3.
\bibitem[CF]. V. P. Camillo, K. R. Fuller, {\it On graded rings with finiteness conditions}, 
Proc. Amer. Math. Soc.
{\bf 86}, no. 1, (1982).
\bibitem[GL]. 
J. Graham, G.I. Lehrer, {\it Cellular algebras}, Inventiones Math. {\bf 123} (1996), 1--34.
\bibitem[EW]. B. Elias, G. Williamson, {\it The Hodge theory of Soergel bimodules}, 
Annals of Mathematics, 
(2) {\bf 180} (2014), no. 3, 1089-1136.
\bibitem[F]. M. Fayers, {\it Decomposition numbers for weight three blocks of symmetric 
groups and Iwahori-Hecke algebras}, Trans. Amer. Math. Soc. {\bf 360} (2008), 1341--1376. 
\bibitem[HuMa1]. J. Hu, A. Mathas, 
{\it Graded cellular bases for the cyclotomic Khovanov-Lauda-Rouquier algebras of type $A$},
Adv. Math., {\bf 225} (2010), 598--642. 
\bibitem[HuMa2]. J. Hu, A. Mathas, 
{\it Seminormal forms and cyclotomic quiver Hecke
algebras of type $A$}, arXiv 1304:0906.
\bibitem[HuMa3]. 
J. Hu, A. Mathas, {\it Graded induction for Specht modules}, Int. Math. Res. Notices {\bf 2012}  (2012)
doi: 10.1093/imrn/rnr058.
\bibitem[Ja].
G. James, {\it The decomposition matrices of $ GL_n(q)$ for $n \leq 10$}, Proc. London Math. Soc., {\bf 60} 
(1990), 225--265.
\bibitem[Ju1]. 
A. A. Jucys, {\it On the Young operators of symmetric groups}, Litousk. Fiz Sb. {\bf 6} (1966)
163--180.
\bibitem[Ju2]. 
A. A. Jucys, {\it Factorisation of Young's projection operators of symmetric groups}, Litousk.
Fiz. Sb. {\bf 11} (1971), l--10.
\bibitem[Ju3].A. A. Jucys, {\it Symmetric polynomials and the centre of the symmetric group ring},
Rep. Mat. Phm. { \bf 5} (1974), 107--112.
\bibitem[KL]. M. Khovanov, A. Lauda, {\it A diagrammatic approach to categorification of quantum groups I}, 
Represent. Theory {\bf 13} (2009), 309--347. 
\bibitem[KM]{Kleshchev-Muth} A. Kleshchev, R. Muth, 
{\it Imaginary Schur-Weyl duality}, arXiv:1312.6104.
\bibitem[KN]{Kleshchev-Nash} A. Kleshchev, D. Nash, 
{\it An interpretation of the Lascoux-Leclerc-Thibon algorithm and graded representation theory}, 
Commun. in Algebra {\bf 38} (2010), 4489--4500.
\bibitem[LLT]{lascoux leclerc thibon} A. Lascoux, B. Leclerc, J.-Y. Thibon,
{\it Hecke algebras at roots of unity and crystal bases of quantum affine 
algebras}, Commun. Math. Phys. {\bf 181} (1996), 205--263.
\bibitem[M81].G. E. Murphy, {\it A New Construction of Young's Seminormal 
Representation of the Symmetric Groups}, Journal of Algebra {\bf 69} (1981), 
287--297.
\bibitem[M83].G. E. Murphy, {\it The Idempotents of the Symmetric Groups and
Nakayama's Conjecture}, Journal of Algebra {\bf 81} (1983), 258--265.
\bibitem[M92].G. E. Murphy, {\it On the Representation Theory of the 
Symmetric
Groups and Associated Hecke Algebras}, Journal of Algebra {\bf 152} (1992),
492--513.
\bibitem[M95].G. E. Murphy, {\it The Representations of Hecke Algebras of type
$ A_n $}, Journal of Algebra {\bf 173} (1995), 97--121. 
\bibitem[Ma]. A. Mathas, 
{\it Iwahori-Hecke Algebras and Schur Algebras of the Symmetric Group},  Univ. Lecture
Series 15, Amer. Math. Soc., 1999.
\bibitem[R]. R. Rouquier, {\it 2-Kac-Moody Algebras}, arXiv:0812.5023.
\bibitem[Ri]. M. J. Richards, {\it Some decomposition numbers for Hecke algebras of general linear
groups}, Math. Proc. of the Camb. Phil. Soc., {\bf 119}, (1996), 383--402.
\bibitem[RH1]. S. Ryom-Hansen, 
Grading the translation functors in type $A$, {\it J. Algebra} {\bf 274} (2004), no. 1, 138--163.
\bibitem[RH2]. S. Ryom-Hansen, 
{\it On the denominators of Young's seminormal basis}, arXiv:0904.4243.	
\bibitem[RH3]. S. Ryom-Hansen, 
{\it Young's seminormal form and simple modules for $ S_n $ in characteristic $p$}, 
Algebras and Representation Theory (2012): 1--23, August 31, 2012.
\bibitem[S97]. M. Sch\"onert et~al.
\newblock { {GAP} --
            {Groups}, {Algorithms}, and {Programming} --
            version 3 release 4 patchlevel 4"}.
\newblock Lehrstuhl D f\"ur Mathematik,
            Rheinisch Westf\"alische
            Technische Hochschule, Aachen, Germany, 1997.
\bibitem[So]. W. Soergel, {\it Kazhdan-Lusztig-Polynome und unzerlegbare Bimoduln \"uber 
Polynomringen}, J. Inst. Math. Jussieu {\bf 6} (2007), no. 3, 501--525.
\bibitem[Wi]. G. Williamson, {\it Schubert calculus and torsion}, 
arXiv:1309.5055v1.
\end{thebibliography}
\end{document}